\documentclass[11pt,letterpaper,reqno]{amsart}
\usepackage{amsmath,amsthm,amsfonts,amssymb,mathtools,algpseudocode}
\usepackage{comment,bookmark,bm,color,caption,subcaption}

\hypersetup{pdfstartview={FitH}}
\def\restrict#1{\raise-.5ex\hbox{\ensuremath|}_{#1}}
\newtheorem{Theorem}{{\bf Theorem}}[section]

\newtheorem{Corollary}[Theorem]{{\bf Corollary}}

\newtheorem{Definition}[Theorem]{{\bf Definition}}
\newtheorem{Lemma}[Theorem]{{\bf Lemma}}

\numberwithin{equation}{section}

\newcommand{\C}{\mathbb{C}}
\newcommand{\ve}{\text{vec}}
\newcommand{\off}{\text{off}}
\newcommand{\Off}{\text{OFF}}
\newcommand{\ubce}{\text{UBCE}}
\newcommand{\diag}{\text{diag}}
\newcommand{\calO}{\mbox{$\mathcal O$}}
\newcommand{\Ol}{\mbox{\Large $\mathcal O$}}
\newcommand{\calR}{\mbox{$\mathcal R$}}
\newcommand{\calRl}{\mbox{\Large $\mathcal R$}}
\newcommand{\calJ}{\mbox{$\mathcal J$}}
\newcommand{\Jl}{\mbox{\Large $\mathcal J$}}
\newcommand{\Jlo}{\mbox{$\Jl_{\!\! \mathcal{O}}$}}
\newcommand{\be}{\mbox{$\mathbf{e}$}}

\begin{document}

\title[Convergence of the complex block Jacobi methods]{Convergence of the complex block Jacobi methods under the generalized serial pivot strategies}
\author{Erna Begovi\'{c}~Kova\v{c}}\thanks{\textsc{Erna Begovi\'{c} Kova\v{c}}, University of Zagreb Faculty of Chemical Engineering and Technology, Maruli\'{c}ev trg 19, 10000 Zagreb, Croatia. \texttt{ebegovic@fkit.hr}}
\author{Vjeran Hari}\thanks{\textsc{Vjeran Hari}, University of Zagreb Faculty of Science, Department of Mathematics, Bijeni\v{c}ka 30, 10000 Zagreb, Croatia. \texttt{hari@math.hr}}

\thanks{E.B.K. is supported in part by Croatian Science Foundation under the project UIP-2019-04-5200.}
\thanks{V.H. is supported in part by Croatian Science Foundation under the project IP-2020-02-2240.}
\date{}

\renewcommand{\subjclassname}{\textup{2020} Mathematics Subject Classification}

\begin{abstract}
The paper considers the convergence of the complex block Jacobi diagonalization methods under the large set of the generalized serial pivot strategies. The global convergence of the block methods for Hermitian, normal and $J$-Hermitian matrices is proven. In order to obtain the convergence results for the block methods that solve other eigenvalue problems, such as the generalized eigenvalue problem, we consider the convergence of a general block iterative process which uses the complex block Jacobi annihilators and operators.
\end{abstract}

\maketitle

\section{Introduction}

The Jacobi method is a well known iterative method for solving the eigenvalue problem (EVP) of the symmetric matrices. Although the Jacobi algorithm is simple, it can be made very efficient~\cite{DrVe08a,DrVe08b}, even on the standard computers. The method is characterized by the asymptotic quadratic and cubic convergence~\cite{Wil62,Hari91a,RheeHa93}. It possesses intrinsic parallelism
which makes it an excellent method for various parallel computers~\cite{Sam71,LukPark89}. One of the most important properties of the Jacobi method is its high relative accuracy on definite~\cite{DeVe92,DoKoMo09} and almost diagonal indefinite~\cite{Matejas09} symmetric matrices. The global convergence of the method is well researched~\cite{FoHe60,Han63,HeZi68,Naz75,Masc95,Hari07,HaBe17}.
Its extension to complex Hermitian matrices is straightforward~\cite{HaBe21}.

The efficiency of the Jacobi method can be improved if the algorithm works on matrix blocks instead of elements. Then, BLAS3 routines can be used to better exploit the cache hierarchy of contemporary computers. Such methods are referred to as the block Jacobi methods.
The cyclic element-wise (or simple) and block Jacobi methods are asymptotically quadratically convergent provided the eigenvalues of the matrix are simple. This contributes to their efficiency. As for the global convergence, many special problems are solved~\cite{YamOksaVajt89,Drmac07,HaSiSi10,HaSiSi14,YamLaKu14}, but still, the global convergence of the element-wise and block method under the general cyclic pivot strategy has not been solved. Therefore, every advancement is valuable, especially if it refers to the block methods. One attempt in this direction~\cite{Hari15} is linked to the global convergence of the real block Jacobi methods under a wider class of cyclic strategies that are weakly equivalent to the serial ones. Later, a much larger class of cyclic strategies was introduced in \cite{HaBe17}. The common name for those strategies is generalized serial strategies. This class encompasses many of the known
cyclic pivot strategies including those used in parallel computing. In \cite{HaBe17} the global convergence of the real element-wise and block Jacobi methods for the symmetric and $J$-symmetric matrices was proved.
In \cite{Hari86} and \cite{HaBe17} it was shown how the theory of Jacobi annihilators and operators can be used to prove convergence of more general real element-wise and block iterative processes.
Here, we generalize the theory from \cite{HaBe17} to work with complex matrices. It now uses the complex Jacobi annihilators and operators which are generalizations of those used in \cite{Hari86} for the element-wise complex Jacobi processes. The obtained results can be used for proving the global convergence of different complex block Jacobi methods for solving ordinary and generalized EVP. We apply them to the block Jacobi methods for solving the EVP of Hermitian, normal and J-Hermitian matrices. Further applications lie in proving the global convergence of the complex block HZ (Hari-Zimmermann)~\cite{Hari19} and CJ (Cholesky-Jacobi)~\cite{Hari18} method. These are complex block Jacobi methods for solving the generalized EVP.

The paper is divided into nine sections. In Section~\ref{sec:method} we consider the complex block Jacobi method for the Hermitian matrix, and in Section~\ref{sec:annihilators} we define and study the block Jacobi annihilators and operators. In Section~\ref{sec:strategies} we describe the set of the generalized serial pivot strategies and prove some results for the block Jacobi operators linked to generalized strategies.
Then we prove the global convergence of the complex block and element-wise Jacobi method for Hermitian matrices under the class of generalized serial strategies. In Section~\ref{sec:normal} we consider the complex block Jacobi method for normal matrices. The Section~\ref{sec:Jacobi-type} is concerned with the general complex block Jacobi-type process. These results are then used in Section~\ref{sec:J-Jacobi} to prove the global convergence of the complex element-wise and block $J$-Jacobi method under the class of generalized strategies. The results of our numerical tests are presented in Section~\ref{sec:num}.
We end the paper with a short conlcusion.

\section{The Complex block Jacobi method}\label{sec:method}

Let $A\in\C^{n\times n}$ be a Hermitian matrix and let
\begin{equation}\label{partition}
\pi=(n_1,n_2,\ldots,n_m), \quad n_1+n_2+\cdots+n_m=n, \quad n_r\geq 1, \quad 1\leq r\leq m,
\end{equation}
be an integer partition of $n$. Then, $\pi$ defines a block-matrix partition $(A_{rs})$ of $A$. Each diagonal block $A_{rr}$, $1\leq r\leq m$, is an $n_r\times n_r$ square, while the other blocks can be rectangular.

The \emph{complex block Jacobi method} applied on $A\in\C^{n\times n}$ is an iterative method,
\begin{equation}\label{jacobistep}
A^{(k+1)}=\left(U^{(k)}\right)^*A^{(k)}U^{(k)}, \quad k\geq0, \quad A^{(0)}=A,
\end{equation}
where the transformation matrices $U^{(k)}$ are unitary \emph{elementary block matrices} of the form
\begin{equation}\label{matrixU}
U^{(k)}=U(i(k),j(k))=\left[
        \begin{array}{ccccc}
          I &  &  &  &  \\
           & U_{ii}^{(k)} &  & U_{ij}^{(k)} &  \\
           &  & I &  &  \\
           & U_{ji}^{(k)} &  & U_{jj}^{(k)} &  \\
           &  &  &  & I \\
        \end{array}
      \right]\begin{array}{c}
                \\
               n_i \\
                \\
               n_j \\
                \\
             \end{array}.
\end{equation}
Indices $i$ and $j$ in~\eqref{matrixU} are called \emph{pivot indices} and they depend on the iteration step $k$. The transformation matrix $U^{(k)}$ of the form~\eqref{matrixU} is determined by its pivot submatrix
\begin{equation}\label{pivsubm}
\widehat{U}_{ij}^{(k)}=\left[\begin{array}{cc}U_{ii}^{(k)}&U_{ij}^{(k)}\\ U_{ji}^{(k)}&U_{jj}^{(k)} \end{array}\right].
\end{equation}
We can write $U^{(k)}=\mathcal{E}(i,j,\widehat{U}_{ij}^{(k)})$, where $\mathcal{E}$ is embedding of $\widehat{U}_{ij}^{(k)}$ into the identity matrix $I_n$.

The goal of the $k$th step of the method~\eqref{jacobistep} is to diagonalize the $(n_i+n_j)\times(n_i+n_j)$ \emph{pivot submatrix} $\widehat{A}_{ij}^{(k)}$ of $A^{(k)}$,
\[
\widehat{A}_{ij}^{(k)}=\left[\begin{array}{cc}
A_{ii}^{(k)} & A_{ij}^{(k)} \\
(A_{ij}^{(k)})^* & A_{jj}^{(k)} \\
\end{array}\right].
\]
For each $k$ the pivot submatrix $\widehat{U}_{ij}^{(k)}$ should satisfy the equation
\begin{equation}\label{pivotstep}
\left(\widehat{U}_{ij}^{(k)}\right)^*\widehat{A}_{ij}^{(k)}\widehat{U}_{ij}^{(k)} =
\left[\begin{array}{cc}\Lambda_{ii}^{(k+1)} & 0 \\ 0 & \Lambda_{jj}^{(k+1)}\end{array}\right],
\end{equation}
for some diagonal matrices $\Lambda_{ii}^{(k+1)}\in\C^{n_i\times n_i}$ and $\Lambda_{jj}^{(k+1)}\in\C^{n_j\times n_j}$. Then $\Lambda_{ii}^{(k+1)}$ and $\Lambda_{jj}^{(k+1)}$
are the diagonal blocks of $A^{(k+1)}$, i.e., we have $A_{ii}^{(k+1)}=\Lambda_{ii}^{(k+1)}$, $A_{jj}^{(k+1)}=\Lambda_{jj}^{(k+1)}$.

Index pair $(i,j)=(i(k),j(k))$ that is used in step $k$ of the method is called \emph{pivot pair}. For the sake of simplicity of notation, when $k$ is implied, it will be omitted. Each block Jacobi method is defined by a \emph{pivot strategy} which determines the way how the pivot pairs are selected.

Denote by $\mathcal{P}_m=\{(r,s)|1\leq r<s\leq m\}$ the set of all pivot positions and let
$\mathbb{N}_0=\{0,1,2,\ldots \}$. The pivot strategy can be described as a function $I:\mathbb{N}_0\rightarrow\mathcal{P}_m$.
If $I$ is a periodic function with period $M=m(m-1)/2$ and
$\{(i(0),j(0)),(i(1),j(1)),\ldots,(i(m-1),j(m-1))\}=\mathcal{P}_m$, then $I$ is a \emph{cyclic} pivot strategy. In other words, a pivot strategy $I$ is cyclic if, during any $M$ successive steps under the strategy $I$, every off-diagonal block is selected and annihilated exactly once.

Let $\Ol(\mathcal{P}_m)$ denote the set of all finite sequences made of the elements of $\mathcal{P}_m$, such that each element of $\mathcal{P}_m$ appears exactly once in each sequence from $\Ol(\mathcal{P}_m)$. Such sequences are called orderings of $\Ol(\mathcal{P}_m)$. Since they
uniquely determine cyclic pivot strategies we will also call them \emph{pivot orderings}.

Each cyclic pivot strategy defines, and is defined by, some pivot ordering $\mathcal{O}\in\Ol(\mathcal{P}_m)$. It has the form
\[
\mathcal{O}=(i(0),j(0)),(i(1),j(1)),\ldots,(i(M-1),j(M-1)).
\]
We can visually depict an ordering $\mathcal{O}\in\Ol(\mathcal{P}_m)$ using the symmetric $m\times m$ matrix $\mathsf{M}_{\mathcal{O}}=(\mathsf{m}_{rs})$ defined by the rule
$$\mathsf{m}_{i(k)j(k)}=\mathsf{m}_{j(k)i(k)}=k, \quad 0\leq k\leq M-1.$$
We can set $\mathsf{m}_{rr}=-1$, for $1\leq r\leq m$.
Since the pairs $(r,r)$ do not appear in $\mathcal{O}$, we can display $*$ instead of $-1$ on the diagonal of $\mathsf{M}$. For illustration, if
\[
\mathcal{O}=(1,2),(3,4),(1,3),(2,4),(1,4),(2,3)\in\Ol(\mathcal{P}_4), \quad \text{then }
\mathsf{M}_{\mathcal{O}}=\left[
                             \begin{array}{cccc}
                               * & 0 & 2 & 4 \\
                               0 & * & 5 & 3 \\
                               2 & 5 & * & 1 \\
                               4 & 3 & 1 & * \\
                             \end{array}
                           \right].
\]

The most common cyclic pivot strategies are the \emph{serial strategies}. They are the \emph{column-cyclic} $I_{\text{col}}=I_{\mathcal{O}_{\text{col}}}$ and the \emph{row-cyclic} strategy $I_{\text{row}}=I_{\mathcal{O}_{\text{row}}}$, where
\begin{align}
\mathcal{O}_{\text{col}} & =(1,2),(1,3),(2,3),(1,4),\ldots,(1,m),(2,m),\ldots,(m-1,m), \label{col}\\
\mathcal{O}_{\text{row}} & =(1,2),(1,3),\ldots,(1,m),(2,3),\ldots,(2,m),\ldots,(m-1,m). \label{row}
\end{align}
In this paper we will consider more general pivot strategies which are referred to as \emph{generalized serial strategies}. They are described in Section~\ref{sec:strategies}.

Note that, once the diagonal blocks $A_{ii}^{(k)}$ and $A_{jj}^{(k)}$ are diagonalized, they will remain diagonal during later iterations. Therefore, it is advisable to apply a preprocessing step where the diagonal blocks of $A$ are diagonalized first, and then the iterative process~\eqref{jacobistep} is applied.

The complex block Jacobi method is said to be \emph{globally convergent} if, for any starting Hermitian block matrix $A$, the sequence $(A^{(k)},k\geq0)$ generated by~\eqref{jacobistep} converges to a fixed diagonal matrix $\Lambda$. We will observe the matrix \emph{off-norm} of the iterates,
\begin{equation}\label{off}
\off(A^{(k)})=\|A^{(k)}-\diag(A^{(k)})\|_F, \quad k\geq0,
\end{equation}
where $\diag(A^{(k)})$ stands for the diagonal part of $A^{(k)}$ and $\|\cdot\|_F$ is the Frobenius norm. The off-norm can be can be understood as the distance from the set of diagonal matrices. Using the block-matrix partition of $A^{(k)}$ we can write the relation~\eqref{off} as
$$\off^2(A^{(k)})=\sum_{\substack{i,j=1 \\ i\neq j }}^m \|A_{ij}^{(k)}\|_F^2+\sum_{i=1}^m \off^2(A_{ii}^{(k)}).$$

In order to ensure the global convergence of the block Jacobi method under the serial pivot strategies, Drma\v{c}~\cite{Drmac07} introduced the class of \emph{UBC} (uniformly bounded cosine) transformations. The matrices from that class are unitary block matrices and the singular values of their diagonal blocks are uniformly bounded by a strictly positive lower bound. The UBC transformations are generalizations of the element-wise (real and complex) Jacobi rotations which satisfy the well-known Forsythe-Henrici condition~\cite{FoHe60}, stating that the cosines of the rotation angles must have a strictly positive uniform lower bound.

The class of UBC matrices was refined in~\cite{HaBe17} using the parameter $\varrho$, $0<\varrho\leq 1$, to get the set $\ubce_{\pi}(\varrho)$ of orthogonal \emph{elementary UBC block matrices} with partition $\pi$. The orthogonal matrix $V=(V_{rs})$ carrying the block-matrix partition defined by $\pi$ is from the class $\ubce_{\pi}(\varrho)$ if it differs from the identity matrix $I_n$ in two
diagonal blocks (say $V_{ii}$, $V_{jj}$) and the corresponding off-diagonal blocks ($V_{ij}$, $V_{ji}$), and satisfies (cf.,~\cite{HaBe17}),
\begin{equation}\label{ubce}
\sigma_{min}(V_{ii})=\sigma_{min}(V_{jj}) \geq \varrho\gamma_{ij}>\varrho\tilde{\gamma}_{n}, \quad (i,j)\in\mathcal{P}_m,
\end{equation}
where
\begin{equation}\label{ubce1}
\gamma_{ij}=\max\left\{\frac{3}{\sqrt{(n_j+1)(4^{n_i}+6n_j-1)}},\left(\begin{array}{c}n_i+n_j \\ n_i\end{array}\right)^{-1/2} \right\}, \quad \tilde{\gamma}_{n} = \frac{3\sqrt{2}}{4^n+26}.
\end{equation}
The first expression within the braces is from \cite{Drmac07} while the second one is obtained from \cite[Lemma~3.3]{Hari21} by setting $\epsilon =0$. The lower bound of $\gamma_{ij}$ involving only the dimension $n$ is easily obtained from the first expression within the braces (see \cite{Hari15}).

The bound~\eqref{ubce}--\eqref{ubce1} was proven for real matrices in~\cite{Hari21}. However, following the lines of the proof of \cite[Lemma~3.3]{Hari21} one finds out that the proof is easily modified to hold for the complex unitary matrices. One has to use the complex Householder matrices, each in the form or the product of an appropriate unitary diagonal matrix and the real Householder reflector. We will consider that problem in Section~4.
Thus, the class $\ubce_{\pi}(\varrho)$ of unitary elementary UBC complex block matrices is defined the same way. One just replaces the word orthogonal with unitary.

\section{Block Jacobi annihilators and operators}~\label{sec:annihilators}

The \emph{Jacobi annihilators and operators} were introduced by Henrici and Zimmermann~\cite{HeZi68} as a tool for proving the global and quadratic convergence of the column-cyclic Jacobi method for symmetric matrices. Later they were used for proving the global convergence of more general Jacobi-type methods~\cite{Hari82,Hari86}, as well as for the block Jacobi methods for symmetric matrices~\cite{Hari15,HaBe17}. In this section we define a class of Jacobi annihilators and operators designed for the block Jacobi methods for the complex Hermitian matrices. They will be referred to as \emph{block Jacobi annihilators and operators}. The definitions given here are from \cite{Hari09}, and they are generalizations of the definitions from~\cite{Hari86} and~\cite{HaBe17}. They are needed here for the sake of completeness.

For an arbitrary matrix $X\in\C^{p\times q}$ we define vectors
\begin{align*}
\textrm{col}(X) & =\left[x_{11} x_{21} \cdots x_{p1} \cdots x_{1q} \cdots x_{pq}\right]^T \in\C^{pq\times1} , \\
\textrm{row}(X) & =\left[ x_{11} x_{12} \cdots x_{1q} \cdots x_{p1} \cdots x_{pq}\right] \in\C^{1\times pq}.
\end{align*}
Let $A\in\C^{n\times n}$ be a block matrix partitioned according to $\pi=(n_1,n_2,\ldots,n_m)$. Set
\[
c_j = \left[
   \begin{array}{c}
     \textrm{col}(A_{1j}) \\
     \textrm{col}(A_{2j}) \\
     \vdots \\
     \textrm{col}(A_{j-1,j}) \\
   \end{array}
 \right], \quad 2\leq j\leq m,
\]
and
\[
r_i=\left[
   \begin{array}{cccc}
     \textrm{row}(A_{i1}) & \textrm{row}(A_{i2}) & \cdots & \textrm{row}(A_{i,i-1}) \\
   \end{array}
 \right], \quad 2\leq i\leq m.
 \]
The column- (row-) vectors $c_j$ ($r_i$) vary in length since they are made of the elements of the
$j$th ($i$th) block-column (block-row) in the strictly upper (lower) triangle of $A$.

Using the vectors $c_j$ and $r_i$, $2\leq i,j\leq m$, we define the function $\ve_{\pi}$, depending on the partition $\pi$,
\begin{equation}\label{NiK}
\ve_{\pi}:\C^{n\times n}\mapsto\C^{2K}, \quad K=N-\sum_{i=1}^m \frac{n_i(n_i-1)}{2}, \quad N=\frac{n(n-1)}{2},
\end{equation}
\begin{equation}\label{vec}
\ve_{\pi}(A) = \left[ c_2^T\ c_3^T\ \cdots c_m^T\ r_2\ r_3\ \cdots r_m\right]^T \in\mathbb{C}^{2K}.
\end{equation}
As we can see, $\ve_{\pi}$ maps the block matrix $A$ to vector $a=\ve_{\pi}(A)$ containing all elements from every off-diagonal block of $A$.
Assuming that $\pi$ is fixed, we can write $\ve$ instead of $\ve_{\pi}$.

The function $\tau$ that maps the position of the block $(i,j)$ in matrix $A$ to its position in the vector $\ve(A)$ is given by
\[
\tau(i,j)=\left\{
              \begin{array}{c}
                (j-1)(j-2)/2+i, \quad \text{for} \ 1\leq i<j\leq m, \\
                \tau(j,i)+M, \quad \text{for} \ 1\leq j<i\leq m. \\
              \end{array}
            \right.
\]
Here, $M=\frac{m(m-1)}{2}$. Obviously, $\tau: \{(i,j)| 1\leq i,j\leq m,\ i\neq j\} \rightarrow \{1,2,\ldots,2M\}$ is a bijection.

Compared to the function $\ve$ in the real case from~\cite{HaBe17}, the function $\ve$ defined by~\eqref{vec} yields a twice longer vector $a=\ve(A)$. This is the case because, for the Hermitian matrices, $\ve(A)$ contains the elements from both the upper and the lower block-triangle of $A$. On the contrary, for symmetric matrices it was enough to take the elements form the blocks in the upper triangle.

The function $\ve$ is a surjection, but not an injection. Therefore, we make use of its restriction to the vector space of complex $n\times n$ matrices carrying the block-matrix partition defined by $\pi$ and with the diagonal blocks equal to zero, denoted by $\mathbb{C}_{0,\pi}^{n\times n}$.
Then
\[
\ve_0= \ve\restrict{\mathbb{C}_{0,\pi}^{n\times n}}, \quad \ve_0:\mathbb{C}_{0,\pi}^{n\times n}\mapsto\mathbb{C}^{2K},
\]
is a bijection.

Moreover, let us define the linear operator
$\mathcal{N}_{ij}:\mathbb{C}^{n\times n}\rightarrow\mathbb{C}^{n\times n}$ as mapping
that sets the blocks on the positions $(i,j)$ and $(j,i)$ of the argument matrix to zero, while the rest of the matrix remains the same.

By employing the functions $\ve$, $\ve_0$, and $\mathcal{N}_{ij}$, we can define \emph{the block Jacobi annihilator} and \emph{the block Jacobi operator} for the complex block matrices. Compared to their analogs for the real case from~\cite{HaBe17}, the main difference is in the sizes of vectors and matrices.

\begin{Definition}\label{def:jan}
Let $\pi =(n_1,\ldots ,n_m)$ be a partition of $n$. Let $\widehat{U}\in\C^{(n_i+n_j)\times(n_i+n_j)}$ be a unitary matrix and let $U=\mathcal{E}(i,j,\widehat{U})$ be the elementary block matrix carrying the matrix block partition defined by $\pi$. The transformation $\calR_{ij}(\widehat{U})$ defined by
\begin{equation}\label{jan}
\calR_{ij}(\widehat{U})(vec(A)) = vec(\mathcal{N}_{ij}(U^*AU)), \quad A\in\C^{n\times n},
\end{equation}
is the complex $ij$-block Jacobi annihilator.

For each pair $(i,j)$, $1\leq i<j\leq m$, the set
\[
\calRl_{ij} =\left\{\calR_{ij}(\widehat{U}) \ \big{|} \ \widehat{U}\in\C^{(n_i+n_j)\times(n_i+n_j)} \text{ unitary} \right\}
\]
is the $ij$-class of the complex block Jacobi annihilators. For any $0<\varrho\leq 1$,
\[
\calRl_{ij}^{\ubce_{\pi}(\varrho)} =\left\{\calR_{ij}(\widehat{U}) \ \big{|} \ \widehat{U}\in\C^{(n_i+n_j)\times(n_i+n_j)} \text{ and } \widehat{U}\in\ubce_{\pi}(\varrho)\right\}
\]
is a subclass of $\calRl_{ij}$.
\end{Definition}

An alternative way of defining the complex $ij$-block Jacobi annihilator would use the relation
\[
\calR_{ij}(\widehat{U})(a) = vec(\mathcal{N}_{ij}(U^*\ve_0^{-1}(a)U)), \quad a\in \C^{2K}, \quad K=N-\sum_{i=1}^m \frac{n_i(n_i-1)}{2}, \quad N=\frac{n(n-1)}{2},
\]
instead of~\eqref{jan}. Actually, $\calR_{ij}(\widehat{U})$ is a $2K\times 2K$ matrix which we also call transformation or mapping because its main purpose is to act on some vector which is connected to some matrix via the relation (\ref{jan}).

A block Jacobi annihilator from $\calRl_{ij}$ can be observed as a generalization of one block Jacobi step. More precisely, a generalization of the mapping that connects two consecutive iterations in the block Jacobi process (\ref{jacobistep}). The difference is the following. In the $k$th iteration of the block Jacobi method, the transformation $U^{(k)}$ is carefully chosen to annihilate the pivot blocks on positions $(i,j)$ and $(j,i)$.
On the other hand, the Jacobi annihilator takes any elementary block matrix $\mathcal{E}(i,j,\widehat{U})$, for an arbitrary unitary $\widehat{U}\in\C^{(n_i+n_j)\times(n_i+n_j)}$, i.e. not necessarily the one annihilating the pivot blocks, then simply sets the pivot blocks to zero.

The structure of a block Jacobi annihilator is given in the following theorem.
\begin{Theorem}\label{tm:jan}
Let $\pi=(n_1,\ldots,n_m)$ be a partition of $n$. Let $(i,j)\in\mathcal{P}_m$ and let $\calR=\calR_{ij}(\widehat{U})\in\calRl_{ij}$ be a complex block Jacobi annihilator. Then $\mathcal{R}$ differs from the identity matrix $I_{2K}$, $K$ from the relation (\ref{NiK}), in exactly $2m+2$ principal submatrices given by the following relations:
$$\mathcal{R}_{\tau(i,j),\tau(i,j)}=0, \quad \mathcal{R}_{\tau(j,i),\tau(j,i)}=0,$$
\begin{align}
\left[
    \begin{array}{cc}
      \mathcal{R}_{\tau(r,i),\tau(r,i)} & \mathcal{R}_{\tau(r,i),\tau(r,j)} \\
      \mathcal{R}_{\tau(r,j),\tau(r,i)} & \mathcal{R}_{\tau(r,j),\tau(r,j)} \\
    \end{array}
  \right] & =\left[
     \begin{array}{cc}
       U_{ii}^T\otimes I_{n_r} & U_{ji}^T\otimes I_{n_r} \\
       U_{ij}^T\otimes I_{n_r} & U_{jj}^T\otimes I_{n_r} \\
     \end{array}
   \right], \quad \text{for} \ 1\leq r\leq m, \label{rel:jantm1}\\
\left[
    \begin{array}{cc}
      \mathcal{R}_{\tau(i,r),\tau(i,r)} & \mathcal{R}_{\tau(i,r),\tau(j,r)} \\
      \mathcal{R}_{\tau(j,r),\tau(i,r)} & \mathcal{R}_{\tau(j,r),\tau(j,r)} \\
    \end{array}
  \right] & =\left[
     \begin{array}{cc}
       I_{n_r}\otimes U_{ii}^* & I_{n_r}\otimes U_{ji}^* \\
       I_{n_r}\otimes U_{ij}^* & I_{n_r}\otimes U_{jj}^* \\
     \end{array}
   \right], \quad \text{for} \ 1\leq r\leq m. \label{rel:jantm2}
\end{align}
Here $\otimes$ denotes the Kronecker product. The matrices from the relations (\ref{rel:jantm1}) and
(\ref{rel:jantm2}) are unitary.
\end{Theorem}

\begin{proof}
Let $a\in \C^{2K}$ be arbitrary and $A=\ve_0^{-1}(a)$. Let us fix the pair $(i,j)$ and let us consider the relation (\ref{jan}). The annihilator $\mathcal{R}$ annihilates the blocks of $a$ in positions $\tau(i,j)$ and $\tau(j,i)$. Therefore, we have $\mathcal{R}_{\tau(i,j),\tau(i,j)}=0$ and $\mathcal{R}_{\tau(j,i),\tau(j,i)}=0$.

Since $A\in \mathbb{C}_{0,\pi}^{n\times n}$, the function $\ve$ acts on $A$ the same way as $\ve_0$.
From the relation (\ref{jan}) we see that $\ve_0^{-1}\left(\calR_{ij}(\widehat{U})(vec(A))\right)$
is also in $\mathbb{C}_{0,\pi}^{n\times n}$. Hence, $A'=\ve_0^{-1}(\calR_{ij}(\widehat{U})(vec(A)))$ is well defined. Note that both $U^*AU$ and
$\mathcal{N}_{ij}(U^*AU)$ are in $\mathbb{C}_{0,\pi}^{n\times n}$. We have
\[
A'= \ve_0^{-1}(\calR_{ij}(\widehat{U})(vec(A))) = \ve_0^{-1}(vec(\mathcal{N}_{ij}(U^*AU)))=
\mathcal{N}_{ij}(U^*AU).
\]
Consequently, for $1\leq r\leq m$, $r\notin\{i,j\}$, we have
\begin{align}
\label{rel:jan1} A'_{ri}=A_{ri}U_{ii}+A_{rj}U_{ji}, & \quad A'_{rj}=A_{ri}U_{ij}+A_{rj}U_{jj}, \\
\label{rel:jan2} A'_{ir}=U_{ii}^*A_{ir}+U_{ji}^*A_{jr}, & \quad A'_{jr}=U_{ij}^*A_{ir}+U_{jj}^*A_{jr}.
\end{align}
From the relation~\eqref{rel:jan2} we obtain
$$A'_{ir}e_l =U_{ii}^*A_{ir}e_l+U_{ji}^*A_{jr}e_l, \quad A'_{jr}e_l =U_{ij}^*A_{ir}e_l+U_{jj}^*A_{jr}e_l, \quad 1\leq l\leq n_r,$$
where $e_l$ stands for the $l$th column of of the identity matrix $I_{n_r}$.
These two equations can be written in the matrix form as
\[
\left[
    \begin{array}{c}
      A_{ir}'e_1 \\[3pt]
      A_{ir}'e_2 \\
      \vdots \\
      A_{ir}'e_{n_r} \\
      A_{jr}'e_1 \\
      A_{jr}'e_2 \\
      \vdots \\
      A_{jr}'e_{n_r} \\
    \end{array}
  \right]=\left[
            \begin{array}{cccc|cccc}
              U_{ii}^* & & &  & U_{ji}^* &   & &  \\
              & U_{ii}^* & &  & &   U_{ji}^* & & \\
              &   & \ddots &    & & & \ddots & \\
              &   &  & U_{ii}^* & & &  & U_{ji}^* \\ \hline
              U_{ij}^* & & &  & U_{jj}^* &   & &  \\
              & U_{ij}^* & &  & &   U_{jj}^* & & \\
              &   & \ddots &   & & & \ddots & \\
              &   &  & U_{ij}^* &  & &  & U_{jj}^* \\
            \end{array}
          \right]\left[
    \begin{array}{c}
      A_{ir}e_1 \\
      A_{ir}e_2 \\[3pt]
      \vdots \\
      A_{ir}e_{n_r} \\ \hline
      A_{jr}e_1 \\
      A_{jr}e_2 \\
      \vdots \\
      A_{jr}e_{n_r} \\
    \end{array}
  \right].
\]
Using another notation, the last relation can be written as
\[
\left[\begin{array}{c} \textrm{col} (A_{ir}')\\ \textrm{col}(A_{jr}')\end{array}\right]
=\left[\begin{array}{cc}
       I_{n_r}\otimes U_{ii}^* & I_{n_r}\otimes U_{ji}^* \\
       I_{n_r}\otimes U_{ij}^* & I_{n_r}\otimes U_{jj}^* \\
     \end{array}\right]
\left[\begin{array}{c} \textrm{col} (A_{ir})\\ \textrm{col} (A_{jr})\end{array}\right].
\]
The blocks $\textrm{col} (A_{ir})$, $\textrm{col} (A_{jr})$ are situated in $a$ at positions
$\tau (i,r)$, $\tau (j,r)$, respectively. Therefore, we have
\[
\left[\begin{array}{cc}
      \calR_{\tau(i,r),\tau(i,r)} & \calR_{\tau(i,r),\tau(j,r)} \\
      \calR_{\tau(j,r),\tau(i,r)} & \calR_{\tau(j,r),\tau(j,r)}
    \end{array}\right]
    =
    \left[\begin{array}{cc}
       I_{n_r}\otimes U_{ii}^* & I_{n_r}\otimes U_{ji}^* \\
       I_{n_r}\otimes U_{ij}^* & I_{n_r}\otimes U_{jj}^* \\
     \end{array}\right]
     =
    \left[\begin{array}{cc}
       I_{n_r}\otimes U_{ii} & I_{n_r}\otimes U_{ij} \\
       I_{n_r}\otimes U_{ji} & I_{n_r}\otimes U_{jj} \\
     \end{array}\right]^*,
\]
and this proves the assertion~\eqref{rel:jantm2}. Using the properties of the Kronecker product
(see \cite[Section~6.5]{FuZe11} it is easy to prove that the matrix from the relation~\eqref{rel:jantm2} is unitary.

Now, consider the relation~\eqref{rel:jan1}. Using again the columns of $I_{n_r}$ we have
\[
e_l^TA'_{ri} =e_l^TA_{ri}U_{ii}+e_l^TA_{rj}U_{ji}, \quad
e_l^TA'_{rj} =e_l^TA_{ri}U_{ij}+e_l^TA_{rj}U_{jj}, \quad 1\leq l\leq n_r,
\]
and after applying Hermitian transpose, we obtain
\[
(A'_{ri})^* e_l =U_{ii}^*A_{ri}^* e_l+U_{ji}^*A_{rj}^* e_l, \quad
(A'_{rj})^* e_l =U_{ij}^*A_{ri}^* e_l+U_{jj}^*A_{rj}^* e_l, \quad 1\leq l\leq n_r.
\]
Thus,
\[
\left[\begin{array}{c}
      (A_{ri}')^*e_1 \\
      (A_{ri}')^*e_2 \\
      \vdots \\
      (A_{ri}')^*e_{n_r} \\
      (A_{rj}')^*e_1 \\
      (A_{rj}')^*e_2 \\
      \vdots \\
      (A_{rj}')^*e_{n_r} \\
    \end{array}
  \right]=\left[\begin{array}{cccc|cccc}
              U_{ii}^* & & &  & U_{ji}^* &   & &  \\
              & U_{ii}^* & &  & &   U_{ji}^* & & \\
              &   & \ddots &    & & & \ddots & \\
              &   &  & U_{ii}^* & & &  & U_{ji}^* \\ \hline
              U_{ij}^* & & &  & U_{jj}^* &   & &  \\
              & U_{ij}^* & &  & &   U_{jj}^* & & \\
              &   & \ddots &   & & & \ddots & \\
              &   &  & U_{ij}^* &  & &  & U_{jj}^* \\
            \end{array}\right]
  \left[\begin{array}{c}
      A_{ri}^*e_1 \\
      A_{ri}^*e_2 \\
      \vdots \\
      A_{ri}^*e_{n_r} \\
      A_{rj}^*e_1 \\
      A_{rj}^*e_2 \\
      \vdots \\
      A_{rj}^*e_{n_r} \\
    \end{array}\right].
\]
Using a shorter notation we obtain
\[
\left[\begin{array}{c} \textrm{col} ((A_{ri}')^*)\\ \textrm{col} ((A_{rj}')^*)\end{array}\right]
=\left[\begin{array}{cc}
       I_{n_r}\otimes U_{ii}^* & I_{n_r}\otimes U_{ji}^* \\
       I_{n_r}\otimes U_{ij}^* & I_{n_r}\otimes U_{jj}^* \\
     \end{array}\right]
\left[\begin{array}{c} \textrm{col} ((A_{ri})^*)\\ \textrm{col} ((A_{rj})^*)\end{array}\right].
\]
This is equivalent to
\begin{eqnarray}
\textrm{col} ((A_{ri}')^*) = & (I_{n_r}\otimes U_{ii}^*)\textrm{col} ((A_{ri})^*)+
(I_{n_r}\otimes U_{ji}^*)\textrm{col} ((A_{rj})^*), \label{rel:janmatrix1} \\
\textrm{col} ((A_{rj}')^*) = & (I_{n_r}\otimes U_{ij}^*)\textrm{col} ((A_{ri})^*)+ (I_{n_r}\otimes U_{jj}^*)\textrm{col} ((A_{rj})^*). \label{rel:janmatrix2}
\end{eqnarray}
Now, we make use of the $n_rn_i\times n_rn_i$ permutation matrix $P$ such that
\[
P^T\textrm{col} ((A_{ri})^*)
=
P^T\left[\begin{array}{c}
      A_{ri}^*e_1 \\
      A_{ri}^*e_2 \\
      \vdots \\
      A_{ri}^*e_{n_r} \\
    \end{array}\right]
    =
    \left[\begin{array}{c}
      \bar{A}_{ri}e_1 \\
      \bar{A}_{ri}e_2 \\
      \vdots \\
      \bar{A}_{ri}e_{n_r} \\
    \end{array}\right] = \textrm{col} (\bar{A}_{ri}),
\]
where $\bar{A}_{ri}$ stands for the complex conjugate of $A_{ri}$.
The permutation $P$ can be described by its columns. If $I_{n_r n_i}=[\be_1,\be_2,\ldots ,\be_{n_r n_i}]$ we have
\begin{eqnarray*}
  P &=& [\be_1\ \be_{n_i+1}\ \be_{2n_i+1}\ \be_{3n_i+1}\ \ldots\ \be_{(n_r-1)n_i+1}\
  \be_2\ \be_{n_i+2}\ \be_{2n_i+2}\ \ldots\ \be_{(n_r-1)n_i+2}\ \ldots  \\
   & & \quad \be_{n_i}\ \be_{2n_i}\ \be_{3n_i}\ \ldots\ \be_{n_r n_i}].
\end{eqnarray*}
From the relation~\eqref{rel:janmatrix1} we have
\begin{eqnarray*}
\textrm{col} ((\bar{A}_{ri}')) &=& P^T\textrm{col} ((A_{ri}')^*) = P^T(I_{n_r}\otimes U_{ii}^*)PP^T\textrm{col} ((A_{ri})^*)+P^T(I_{n_r}\otimes U_{ji}^*)\,PP^T\textrm{col} ((A_{rj})^*)\\
&=&  P^T(I_{n_r}\otimes U_{ii}^*)P\ \textrm{col} (\bar{A}_{ri})+P^T(I_{n_r}\otimes U_{ji}^*)P\ \textrm{col} (\bar{A}_{rj}).
\end{eqnarray*}
Applying the complex conjugate operator to the latest equation we conclude that
\begin{align}
\mathcal{R}_{\tau(r,i),\tau(r,i)} =P^T(I_{n_r}\otimes U_{ii}^T)P
& =P^T\left[\begin{array}{ccc}
              U_{ii}^T &  &  \\
               & \ddots &  \\
               &  & U_{ii}^T  \\
            \end{array}\right]P
= U_{ii}^T\otimes I_{n_r}, \label{rel:janmatrix3} \\
\mathcal{R}_{\tau(r,i),\tau(r,j)} =P^T(I_{n_r}\otimes U_{ji}^T)P
& =P^T\left[\begin{array}{ccc}
              U_{ji}^T &  &  \\
               & \ddots &  \\
               &  & U_{ji}^T  \\
            \end{array}\right]P
=
U_{ji}^T\otimes I_{n_r}. \label{rel:janmatrix4}
\end{align}
The right-most equations in the relations (\ref{rel:janmatrix3}) and (\ref{rel:janmatrix4}) follow after carefully inspecting how the similarity transformation with $P$ changes the elements of the matrices $I_{n_r}\otimes U_{ii}^T$ and $I_{n_r}\otimes U_{ji}^T$.
In the same way, from the relation~\eqref{rel:janmatrix2}, we obtain
\begin{equation}\label{rel:janmatrix56}
\mathcal{R}_{\tau(j,r),\tau(i,r)} =U_{ij}^T\otimes I_{n_r}, \quad \mathcal{R}_{\tau(j,r),\tau(j,r)} =U_{jj}^T\otimes I_{n_r}.
\end{equation}
The relations~\eqref{rel:janmatrix3}--\eqref{rel:janmatrix56} yield the assertion~\eqref{rel:jantm1}. To see that the matrix from~\eqref{rel:jantm1} is unitary, one uses the basic properties of the Kronecker product and the fact that for any unitary matrix $W$ the matrices
$\bar{W}$ and $W^T$ are also unitary.

All the other blocks of $\mathcal{R}$, except those considered in this proof, are as in the identity matrix $I_{2K}$.
\end{proof}

From the Theorem~\ref{tm:jan} we conclude that each complex block Jacobi annihilator $\calR_{ij}(\widehat{U})$ is, up to a permutational similarity transformation, a direct sum of the null matrix of dimension $2n_in_j$, the identity matrix of dimension $2K-2n_in_j-2mn(n_i+n_j)$ and some unitary matrix of order $2mn(n_i+n_j)$. However, if $\pi=(n_i,n_j)$, then $\calR_{ij}(\widehat{U})$ is just the null matrix of dimension $2n_in_j$.

The complex block Jacobi annihilators are the building blocks of the complex block Jacobi operators. In the same way as the block annihilators generalize one step of the complex block Jacobi method, the block operators generalize one cycle of the block method.

\begin{Definition}\label{def:jop}
Let $\pi=(n_1,\ldots,n_m)$ be a partition of $n$ and let
\[
\mathcal{O} = (i_0,j_0),(i_1,j_1),\ldots,(i_{M-1},j_{M-1})\in\Ol(\mathcal{P}_m), \quad M=\frac{m(m-1)}{2}.
\]
Then
\[
\Jlo = \{\calJ \ \big{|} \ \calJ=\calR_{i_{M-1}j_{M-1}}\ldots\calR_{i_1j_1}\calR_{i_0j_0},
\ \calR_{i_kj_k}\in\calRl_{i_k j_k}, \ 0\leq k\leq M-1\}
\]
is the \emph{class of the complex block Jacobi operators} associated with the pivot ordering $\mathcal{O}$.

The $2K\times2K$ matrices $\calJ\in\Jlo$, $K=N-\sum_{i=1}^m \frac{n_i(n_i-1)}{2}$, $N=\frac{n(n-1)}{2}$, are called the complex block Jacobi operators.
If $\calR_{i_kj_k}\in\calRl_{i_kj_k}^{\ubce_{\pi}(\varrho)}$ holds for all $0\leq k\leq M-1$ then we write $\Jlo^{\ubce_{\pi}(\varrho)}$ instead of $\Jlo$.
\end{Definition}

\section{Convergence of the complex block Jacobi method under the generalized serial pivot strategies}\label{sec:strategies}

We consider convergence of the complex block Jacobi method under the large class of the generalized serial pivot strategies. This section heavily depends on the theory developed in \cite{HaBe17} and \cite{BegovicPhD}. To keep the exposition as brief as possible we will frequently cite the results from \cite{HaBe17}. We note that definition of the pivot ordering is independent of the \emph{core algorithm} which is used to compute the pivot submatrix $\widehat{U}_k$ of the transformation matrix at each step $k$, $k\geq 0$.

Before we define the set of the generalized serial pivot orderings, which determine the generalized serial pivot strategies, we need to establish certain relations between the pivot orderings.

An \emph{admissible transposition} in a pivot ordering $\mathcal{O}$ is any transposition of two adjacent terms
\[
(i_r,j_r),(i_{r+1},j_{r+1})\rightarrow(i_{r+1},j_{r+1}),(i_r,j_r),
\]
provided that the sets $\{i_r,j_r\}$ and $\{i_{r+1},j_{r+1}\}$ are disjoint.
Two sequences $\mathcal{O},\mathcal{O}'\in\mathcal{\Ol}(\mathcal{S})$, $\mathcal{S}\subseteq \mathcal{P}_n$ are \emph{equivalent} if one can be obtained from the other by a finite set of admissible transpositions. Then we write $\mathcal{O}\sim\mathcal{O}'$.
Further on, $\mathcal{O}$ and $\mathcal{O}'$ are
\begin{itemize}
\item \emph{shift-equivalent} ($\mathcal{O}\stackrel{\mathsf{s}}{\sim}\mathcal{O}'$) if $\mathcal{O}=[\mathcal{O}_1,\mathcal{O}_2]$ and $\mathcal{O}'=[\mathcal{O}_2,\mathcal{O}_1]$, where $[ \ , \  ]$ stands for concatenation of the sequences,
\item \emph{weak equivalent} ($\mathcal{O}\stackrel{\mathsf{w}}{\sim}\mathcal{O}'$) if there exist a positive integer $t$, and sequences $\mathcal{O}_i\in\mathcal{\Ol}(\mathcal{S})$, $0\leq i\leq t$, such that every two adjacent terms in the sequence $\mathcal{O}=\mathcal{O}_0,\mathcal{O}_1,\ldots,\mathcal{O}_t=\mathcal{O}'$ are equivalent or shift-equivalent,
\item \emph{permutation equivalent} ($\mathcal{O}\stackrel{\mathsf{p}}{\sim}\mathcal{O}'$ or $\mathcal{O}'=\mathcal{O}(\mathsf{q})$) if there is a permutation $\mathsf{q}$ such that
    $\mathcal{O}=(i_0,j_0),(i_1,j_1),\ldots, (i_r,j_r)$,
    $\mathcal{O}'=(\mathsf{q}(i_0),\mathsf{q}(j_0)),(\mathsf{q}(i_1),\mathsf{q}(j_1)),\ldots,
    (\mathsf{q}(i_r),\mathsf{q}(j_r))$, where $r+1$ is number of pairs in $\mathcal{O}$,
\item \emph{reverse} ($\mathcal{O}'=\mathcal{O}^{\leftarrow}$) if $\mathcal{O}=(i_0,j_0),(i_1,j_1),\ldots,(i_{r},j_{r})$ and $\mathcal{O}'=(i_{r},j_{r}),\ldots,(i_1,j_1),(i_0,j_0)$.
\end{itemize}
It is easy to check that $\sim$, $\stackrel{\mathsf{s}}{\sim}$, $\stackrel{\mathsf{w}}{\sim}$, $\stackrel{\mathsf{p}}{\sim}$ are reflexive, symmetric and transitive, hence they are equivalence relations on the set $\mathcal{\Ol}(\mathcal{S})$. Note that permutation $\mathsf{q}$ changes the partition $\pi$ and consequently the block-matrix partition.

\emph{Serial orderings with permutations} from~\cite{BegovicPhD,HaBe17} are the orderings of $\mathcal{\Ol}(\mathcal{P}_m)$ where the pivot pairs are arranged (as positions in the matrix) column-to-column or row-to-row, like it is the case with the serial ordering~\eqref{col}, \eqref{row}. However, inside each column (row) the pivot positions come in an arbitrary ordering. These orderings are formally defined as follows.

Denote the set of all permutations of the set $\{l_1,l_1+1,\ldots,l_2\}$ by $\Pi^{(l_1,l_2)}$, and set
\begin{align*}
\mathcal{B}_c^{(m)} = \big{\{} \mathcal{O}\in\mathcal{\Ol}(\mathcal{P}_m) \ \big{|} & \ \mathcal{O}= (1,2),(\tau_{3}(1),3),(\tau_{3}(2),3),\ldots,(\tau_{m}(1),m),\ldots, (\tau_{m}(m-1),m), \\
&  \quad \tau_{j}\in\Pi^{(1,j-1)}, \ 3\leq j\leq m \big{\}},
\end{align*}
\begin{align*}
\mathcal{B}_r^{(m)} = \big{\{}\mathcal{O}\in\mathcal{\Ol}(\mathcal{P}_m) \ \big{|} & \ \mathcal{O}=(m-1,m),(m-2,\tau_{m-2}(m-1)),(m-2,\tau_{m-2}(m)), \\
& \quad \ldots,(1,\tau_{1}(2)),\ldots,(1,\tau_{1}(m)), \quad \tau_{i}\in\Pi^{(i+1,m)}, \ 1\leq i\leq m-2 \big{\}}.
\end{align*}
Then, the set of the serial pivot orderings with permutations is given by
\[
\mathcal{B}_{sp}^{(m)}= \mathcal{B}_c^{(m)} \cup \overleftarrow{\mathcal{B}}_c^{(m)} \cup \mathcal{B}_r^{(m)} \cup \overleftarrow{\mathcal{B}}_r^{(m)},
\]
where
\[
\overleftarrow{\mathcal{B}}_c^{(m)}=\{\overleftarrow{\mathcal{O}},\mathcal{O}\in \mathcal{B}_c^{(m)}\}, \quad \overleftarrow{\mathcal{B}}_r^{(m)}=\{\overleftarrow{\mathcal{O}},\mathcal{O}\in \mathcal{B}_r^{(m)}\}.
\]
The class $\mathcal{B}_{sp}^{(m)}$ can be enlarged by adding all orderings $\mathcal{O}\in\Ol(\mathcal{P}_m)$ that are related to some orderings from $\mathcal{B}_{sp}^{(m)}$ with one or more equivalence relations. We say that two orderings $\mathcal{O},\mathcal{O}'\in\mathcal{\Ol}(\mathcal{P}_m)$ are connected by a chain of equivalence relations if there exist a sequence of orderings
$\mathcal{O}=\mathcal{O}_1,\ldots ,\mathcal{O}_r=\mathcal{O}'\in \mathcal{\Ol}(\mathcal{P}_m)$, $r\geq 2$, such that each two neighboring terms are linked with one of the relations
$\sim$, $\stackrel{\mathsf{s}}{\sim}$, $\stackrel{\mathsf{w}}{\sim}$ or $\stackrel{\mathsf{p}}{\sim}$. Fortunately, every such chain of equivalence relations can be reduced
to the \emph{canonical form} (see \cite{HaBe17}).
The chain connecting $\mathcal{O}$ and $\mathcal{O}''\in\mathcal{\Ol}(\mathcal{P}_m)$ is in the canonical form if it looks like
\[
\mathcal{O}\stackrel{\mathsf{p}}{\sim}\mathcal{O}'\stackrel{\mathsf{w}}{\sim}\mathcal{O}'' \quad \text{or} \quad \mathcal{O}\stackrel{\mathsf{w}}{\sim}\mathcal{O}'\stackrel{\mathsf{p}}{\sim}\mathcal{O}''.
\]
Note that by choosing the permutation to be the identity we obtain $\mathcal{O}\stackrel{\mathsf{w}}{\sim}\mathcal{O}''$ and by choosing $\mathcal{O}'=\mathcal{O}''$ or $\mathcal{O}=\mathcal{O}'$ we obtain $\mathcal{O}\stackrel{\mathsf{p}}{\sim}\mathcal{O}''$.

The set of the \emph{generalized serial pivot orderings} is defined as
$$\mathcal{B}_{sg}^{(m)} = \big{\{}\mathcal{O}\in\Ol(\mathcal{P}_m) \ \big{|} \ \mathcal{O}\stackrel{\mathsf{p}}{\sim}\mathcal{O'}\stackrel{\mathsf{w}}{\sim}\mathcal{O''} \ \text{or} \ \mathcal{O}\stackrel{\mathsf{w}}{\sim}\mathcal{O'}\stackrel{\mathsf{p}}{\sim}\mathcal{O''}, \ \mathcal{O''}\in\mathcal{B}_{sp}^{(m)}\big{\}}.$$
This set of pivot orderings is very broad and contains the ordering that are noticeably different from the serial pivot orderings with permutations. For example, from the set $\mathcal{B}_{sg}^{(5)}$ we have the ordering
\[
\mathcal{O}=(1,4),(4,5),(1,3),(2,4),(3,5),(2,3),(1,5),(1,2),(3,4),(2,5)\in\mathcal{B}_{sg}^{(5)}
\]
with matrix representation
\[
\mathsf{M}_{\mathcal{O}}=
\left[
  \begin{array}{ccccc}
    * & 7 & 2 & 0 & 6 \\
    7 & * & 5 & 3 & 9 \\
    2 & 5 & * & 8 & 4 \\
    0 & 3 & 8 & * & 1 \\
    6 & 9 & 4 & 1 & * \\
  \end{array}
\right].
\]

\subsection{The convergence results}

Now, we can consider the global convergence of the complex block Jacobi method under the generalized serial pivot strategies. Theorem~\ref{tm:sg_jop} refers to the process defined by the complex block Jacobi annihilators. It is then used in the proof of Theorem~\ref{tm:sg} which deals with the complex block Jacobi method.

\begin{Theorem}\label{tm:sg_jop}
Let $\pi=(n_1,\ldots,n_m)$ be a partition of $n$ and $\mathcal{O}\in\mathcal{B}_{sg}^{(m)}$.
Suppose that the chain connecting $\mathcal{O}$ and $\mathcal{O''}\in\mathcal{B}_{sp}^{(m)}$ in the canonical form contains $d$ shift equivalences. Let $0<\varrho\leq1$ and let $\mathcal{J}_1$, $\mathcal{J}_2$,\ldots,$\mathcal{J}_{d+1}$ be any $d+1$ block Jacobi annihilators from $\Jlo^{\ubce_{\pi}(\varrho)}$. Then, there exist constants $\mu_{\pi,\varrho}$ depending only on $\pi$ and $\varrho$, and $\widetilde{\mu}_{n,\varrho}$ depending only on $n$ and $\varrho$, such that
$$\|\mathcal{J}_{d+1}\cdots\mathcal{J}_2\mathcal{J}_1\|_2 \leq \mu_{\pi,\varrho}, \quad 0\leq\mu_{\pi,\varrho}<\widetilde{\mu}_{n,\varrho}<1.$$
\end{Theorem}

\begin{proof}
The proof follows the lines of the proof of Theorem 3.9 from~\cite{HaBe17}. The only difference comes from the fact that block Jacobi annihilators are complex, they are defined by Definition~\ref{def:jan}. The proof also makes use of Theorem~\ref{tm:jan} and relations~\eqref{ubce}, \eqref{ubce1}.
\end{proof}

In Theorem~\ref{tm:sg_jop} we used the spectral matrix norm. The same notation $\|\cdot\|_2$ is used for the Euclidian vector norm.

\begin{Theorem}\label{tm:sg}
Let $\pi=(n_1,\ldots,n_m)$ be a partition of $n$ and $\mathcal{O}\in\mathcal{B}_{sg}^{(m)}$.
Suppose that the chain connecting $\mathcal{O}$ and $\mathcal{O''}\in\mathcal{B}_{sp}^{(m)}$ in the canonical form contains $d$ shift equivalences. Let $A\in\mathbb{C}^{n\times n}$ be a Hermitian matrix carrying the block-matrix partition defined by $\pi$. Let $A'$ be obtained from $A$ by applying $d+1$ sweeps of the cyclic block Jacobi method defined by the pivot strategy $I_{\mathcal{O}}$. If all transformation matrices are from the class $\ubce_{\pi}(\varrho )$, $0<\varrho\leq 1$, then there exist constants $\mu_{\pi,\varrho}$ depending only on $\pi$ and $\varrho$, and $\widetilde{\mu}_{n,\varrho}$ depending only on $n$ and $\varrho$, such that
$$\off(A')\leq \mu_{\pi,\varrho} \off(A), \quad 0\leq \mu_{\pi,\varrho}<\widetilde{\mu}_{n,\varrho}<1.$$
\end{Theorem}

\begin{proof}
Let $\calO=(i_0,j_0),(i_1,j_1),\ldots,(i_{M-1})\in\mathcal{B}_{sg}^{(m)}$. Let us observe how the elements of the off-diagonal blocks are changing. Using the complex block Jacobi annihilators we can write the $k$th step of the complex block Jacobi method as
\[
a^{(k+1)}=\calR_{i_kj_k} a^{(k)}, \quad \calR_{i_kj_k}=\calR_{i_kj_k}(\widehat{U}_{i(k),j(k)}^{(k)}), \quad \calR_{i_kj_k}\in\calRl_{i_kj_k}^{\ubce_{\pi}(\varrho)},
\]
where $i_k=i(k)$, $j_k=j(k)$, $a^{(k)}=\ve(A^{(k)})$, $a^{(k+1)}=\ve(A^{(k+1)})$, $k\geq 0$. In the same way we can observe the $t$th cycle as
$$a^{((t+1)N)}=\calJ_{t+1} a^{(tN)},$$
$$\calJ_{t+1}= \calR_{i_{(t+1)N-1}j_{(t+1)N-1}}\cdots \calR_{i_{tN+1}j_{tN+1}}\calR_{i_{tN}j_{tN}}\in\Jlo^{\ubce_{\pi}(\varrho)}, \ t\geq 0.$$
After $d+1$ cycles we have
\begin{equation}\label{rel:sg}
a'=\mathcal{J}_{d+1}\cdots\mathcal{J}_2\mathcal{J}_1 a, \quad \calJ_1,\calJ_2,\ldots,\calJ_{d+1}\in\Jlo^{\ubce_{\pi}(\varrho)}.
\end{equation}

Applying Theorem~\ref{tm:sg_jop} to the equation~\eqref{rel:sg} we obtain
\begin{equation*}
\|a'\|_2  \leq \|\mathcal{J}_{d+1}\cdots\mathcal{J}_2\mathcal{J}_1\|_2\|a\|_2
 \leq \mu_{\pi,\varrho}\|a\|_2, \quad 0\leq \mu_{\pi,\varrho}<\widetilde{\mu}_{n,\varrho}<1. \label{rel:sgin}
\end{equation*}
This completes the proof since $\off(A)=\|a\|_2$, $\off(A')=\|a'\|_2$.
\end{proof}

Now, we can prove the global convergence of the complex block Jacobi method for Hermitian matrices.
The proof uses ideas from \cite{HaBe17} and \cite{{Hari21}}. We will first consider the most common case when the eigenvalues of the matrix are simple. If the eigenvalues are multiple, we will assume that ultimately the diagonal elements converging to a multiple eigenvalue are contained within one diagonal block. We assume that the core algorithm uses any convergent complex element-wise Jacobi method. In particular, the cyclic one using the generalized serial strategy (see Corollary~\ref{cor:glconv} below), or an optimal one, say the one under the classical pivot strategy. These methods are very good choice because element-wise Jacobi methods are very accurate and efficient on nearly diagonal Hermitian matrices. Note that during the iteration, most of the time the pivot submatrices are nearly diagonal.

The general case, which includes the case where the other eigenvalue methods serve as the core algorithm, and the case where multiplicities of the eigenvalues are arbitrary, is more complicated. We will not
delve into this problem here, but we are sure the proof is similar to the one presented in~\cite{{Hari21}} for the real block HZ method.

Let us consider how the elementary UBC block matrix is obtained from the eigenvector matrix of the pivot submatrix $\widehat{A}_{ij}^{(k)}$. We consider the $k$th step of the complex block Jacobi method. We omit the superscript $(k)$ and assume $i=i(k)$, $j=j(k)$.
Let $\widehat{U}_{ij}$ be as in the relations (\ref{pivsubm}) and (\ref{pivotstep}),
\[
\widehat{U}_{ij}=\left[\begin{array}{cc}U_{ii}&U_{ij}\\ U_{ji}&U_{jj} \end{array}\right].
\]
To obtain the UBC unitary matrix from $\widehat{U}_{ij}$ one applies the QR algorithm with column pivoting to the $n_i\times (n_i+n_j)$ matrix $[U_{ii}\ U_{ij}]$ (see~\cite{Drmac07}). To this end one can use either the complex rotations or the complex Householder matrices. Let us describe the process using the latter choice.
Each complex Householder matrix $H$ can be written as a product of an appropriate unitary diagonal matrix $\Phi$ and the real Householder reflector $\mathsf{H}$, $H=\Phi\mathsf{H}$. The process uses $n_i$ steps and can be described by the equation
$$\left[\begin{array}{cc} I_{n_i-1} & \\ & H_{n_i}^*\end{array}\right]    \cdots   \left[\begin{array}{cc} 1 & \\ & H_2^*\end{array}\right]
        H_1^*[U_{ii}\ U_{ij}]I_{11'}I_{22'} \cdots I_{n_i n_i'}.$$
Here, $H_r$ is the Householder matrix of order $n_i-r+1$, $\mathsf{H}_{n_i}=[-1]$, and
$I_{rr'}$ is the transposition matrix of order $n_i+n_j$ which is used to swap the columns $r$ and $r'$. Thus, $1\leq r'\leq n_i+n_j$. In \cite{Drmac07} it has been shown that the unitary matrix
\[
\widetilde{\widehat{U}}_{ij} =\widehat{U}_{ij} \widehat{P}_{ij}, \quad \widehat{P}_{ij}=I_{11'}I_{22'} \cdots I_{n_i n_i'}
\]
satisfies~\eqref{ubce} as discussed earlier, therefore, it is a very good choice. The proof for the complex case is identical. In our convergence analysis we will use modification of the above algorithm for computing $\widehat{P}_{ij}$. It is based on the following definition (which is the same as  \cite[Definition~3.2]{Hari21}).

\begin{Definition}\label{attr_R}
The core algorithm has property $\mathsf{R}$ if it replaces the permutation matrix $\widehat{P}_{ij}=I_{11'}I_{22'}\cdots I_{n_in_i'}$ with  $\widetilde{\widehat{P}}_{ij}$ using the following procedure: if $\tilde{r}$ is the smallest integer from the set $\{1,2,\ldots ,n_i\}$ such that $\tilde{r}'>n_i$ then set
$\widetilde{\widehat{P}}_{ij} = I_{\tilde{r}\tilde{r}'}\cdots I_{n_in_i'}$. If such $\tilde{r}$ does not exist, set $\widetilde{\widehat{P}}_{ij} = I_{n_{i}+n_{j}}$.
\end{Definition}

In particular, $\widetilde{\widehat{P}}_{ij}$ is the identity matrix provided that
$\widehat{P}_{ij}=\mbox{diag}(\mathsf{P}_{ij},I_{n_j})$ for some permutation matrix $\mathsf{P}_{ij}$ of order $n_i$. We will also say that the core algorithm has attribute $\mathsf{R}$.

\begin{Theorem}\label{tm:glconv}
Let $A$ be a Hermitian matrix and $(A^{(k)}, k\geq 0)$ be the sequence obtained
by applying the complex block Jacobi method to $A$ under any generalized serial pivot strategy.
Let the core algorithm be any convergent complex element-wise Jacobi method.
Let $0<\varrho\leq 1$ and let the transformations matrices $U^{(k)}$ be from the class $\ubce_{\pi}(\varrho )$. Then the assertions (i)--(iii) hold.
\begin{itemize}
\item[(i)] $\lim_{k\rightarrow\infty}\off (A^{(k)})=0$, i.e.\ the method converges to diagonal form,
\item[(ii)] If the eigenvalues of $A$ are simple and the core algorithm has attribute $\mathsf{R}$, then $\lim_{k\rightarrow\infty} A^{(k)}$ exists, and it is the diagonal matrix of the eigenvalues of $A$,
\item[(iii)] If the core algorithm has attribute $\mathsf{R}$ and if
  there is $k_0\geq 0$ such that, for $k\geq k_0$, all diagonal elements affiliated with a multiple eigenvalue are contained within the same diagonal block, then $\lim_{k\rightarrow\infty} A^{(k)}$ exists and it is the diagonal matrix of the eigenvalues of $A$.
\end{itemize}
\end{Theorem}

\begin{proof}
\begin{itemize}
\item[(i)] The sequence $(\off(A^{(k)}), k\geq 0)$ is nondecreasing and therefore convergent. From Theorem~\ref{tm:sg} we conclude that it contains a subsequence converging to zero. Hence, $(\off(A^{(k)}), k\geq 0)$ converges to zero.
\item[(ii)] Since the eigenvalues of $A$ are simple, for large enough $k$ the transformation matrices computed by the Jacobi method are nearly diagonal. This means that the matrix $\widehat{P}_{ij}$, $i=i(k)$, $j=j(k)$, will have the form $\widehat{P}_{ij}=\mbox{diag}(\mathsf{P}_{n_i},I_{n_j})$.
Therefore, $\widetilde{\widehat{P}}_{ij}=I_{n_i+n_j}$. In other words, for any $0<\varrho\leq 1$ there is a $k_0$ such that for $k\geq k_0$ the core algorithm delivers the transformation matrices that are from the class $\ubce_{\pi}(\varrho )$.

The rest of the proof can use the known formulas for the element-wise complex Jacobi method
\begin{eqnarray}
\lefteqn{\hspace{4ex}
\widehat{U}_{ij} = \widehat{R}_{ij}=
\left[\begin{array}{cc}
 \cos\phi_{ij} &   -e^{\imath\alpha_{ij}}\sin\phi_{ij}  \\
 e^{-\imath\alpha_{ij}}\sin\phi_{ij} &  \cos\phi_{ij}
 \end{array}\right], }\label{compl_rot_mat}
\\
&&\alpha_{ij} =\arg (a_{ij}),\quad \tan 2\phi_{ij} =\frac{2|a_{ij}|}{a_{ii}-a_{jj}},\quad -\frac{\pi}{4}\leq \phi_{ij} \leq \frac{\pi}{4}, \label{angles}\\
&&a_{ii}'=a_{ii}+\tan\phi_{ij} |a_{ij}|, \quad a_{jj}'=a_{jj}-\tan\phi_{ij} |a_{ij}|. \label{simple_formulas}
\end{eqnarray}
The formulas~\eqref{simple_formulas} imply that, ultimately, affiliation of the diagonal elements to the eigenvalues is preserved. Namely, if no step of the core algorithm can change the affiliation then also the block step cannot change the affiliation of the diagonal elements.

Another argument for the same claim is actually given above: if the unitary transformation matrix $U^{(k)}$ is sufficiently close to a unitary diagonal matrix then it does not change affiliation of the diagonal elements of $A^{(k)}$ to the eigenvalues of $A$.

\item[(iii)] The proof is similar to the proof of (ii). The difference comes from the fact that now
$(1,2)$ and $(2,1)$ blocks of each $\widehat{U}_{ij}$ tend to zero as $k$ increases. This is implied by the fact that, ultimately, the diagonal elements of the $i$th and $j$th diagonal block are affiliated to different eigenvalues. So, for large enough $k$, we have $\widehat{P}_{ij}=\mbox{diag}(\mathsf{P}_{ij},I_{n_j})$ for some permutation $\mathsf{P}_{ij}$ of order $n_i$, and consequently we will have $\widetilde{\widehat{P}}_{ij}=I_{n_i+n_j}$.

Note that for large $k$, $A_{ii}$ and $A_{jj}$ are nearly diagonal.
If some diagonal elements of $A_{ii}$ are close to a multiple eigenvalue, one can show, using permutations, that $U_{ii}$ has to be such that $U_{ii}^*A_{ii}U_{ii}$ and  $A_{ii}$ are close to each other. The same will hold for $A_{jj}$.
Thus, for $k$ large enough, the affiliation of the diagonal elements will be preserved.
\end{itemize}
\end{proof}

\begin{Corollary}\label{cor:glconv}
The element-wise complex Jacobi method for solving the eigenvalue problem for Hermitian matrices is
globally convergent provided that it is defined by the generalized serial pivot strategy.
\end{Corollary}

\begin{proof}
Choosing $\pi =(1,1,\ldots ,1)$, one obtains that the transformation matrices are complex rotations
defined by the relations~\eqref{compl_rot_mat}--\eqref{simple_formulas}. The relations~\eqref{ubce1} and~\eqref{ubce} hold with $\gamma_{ij}=\sqrt{2}/2$, which is in accordance with the relation~\eqref{angles}. Thus also implies $\cos\phi_{ij}\geq \sqrt{2}/2 \geq \varrho\sqrt{2}/2$, for any $0<\varrho\leq 1$.
From Theorem~\ref{tm:sg} we obtain $\lim_{k\rightarrow\infty}\off(A^{(k)})=0$. Then, one uses the formulas (\ref{simple_formulas}) to show that the diagonal elements converge to the eigenvalues.
\end{proof}

\section{An application to normal matrices}\label{sec:normal}

In this section, we recall that the block Jacobi method for solving the EVP of Hermitian matrices can be used to solve the eigenvalue problem of a normal matrix. We provide here some convergence results.

Matrix $A$ is normal if it satisfies
\[
A^*A=AA^*.
\]
An extensive list of interesting properties of the normal matrix can be found in~\cite{GJSW87} and~\cite{ElIk98}. In particular, the normal matrix $A$ is unitarily diagonalizable. That is, there is a unitary matrix $W$ such that $W^*AW=\Lambda$, where $\Lambda$ is a diagonal matrix.

Let $A\in\C^{n\times n}$ be a block matrix defined by the partition $\pi=(n_1,\ldots,n_m)$ from~\eqref{partition}. We can write $A$ as
\begin{equation}\label{normalA}
A=B+\imath C,
\end{equation}
\begin{equation}\label{normalBC}
B=\frac{1}{2}(A+A^*), \quad C=\frac{1}{2}(A^*-A)\imath.
\end{equation}
The matrices $B$ and $C$ from~\eqref{normalBC} are Hermitian. The relation (\ref{normalA}) is the unique representation of a square matrix $A$ into the sum of a Hermitian and skew-Hermitian matrix. For this reason $B$ ($\imath C$) is called the Hermitian (skew-Hermitian) part of $A$.

The matrices $B$ and $C$ carry the same block-matrix partition as $A$. It is easy to check that $A$, $B$, and $C$ commute. Since all three commuting matrices $A$, $B$, $C$ are normal, they are simultaneously unitarily diagonalizable. Therefore, there is a unitary matrix $U$ such that
\begin{equation}\label{ABCdiag}
U^*AU=\Lambda, \quad U^*BU=\Lambda_B, \quad U^*CU=\Lambda_C,
\end{equation}
where
\[
\Lambda=\text{diag}(\lambda_1,\lambda_2,\ldots,\lambda_n), \quad \Lambda_B=\text{diag}(\beta_1,\beta_2,\ldots,\beta_n), \quad \Lambda_C=\text{diag}(\gamma_1,\gamma_2,\ldots,\gamma_n).
\]
Since $B$ and $C$ are Hermitian such are $\Lambda_B$ and $\Lambda_C$. It means that the diagonal entries of $\Lambda_B$ and $\Lambda_C$ are real, and consequently, $\Lambda_B$ and $\Lambda_C$ are real diagonal matrices. The relations~\eqref{normalA} and \eqref{ABCdiag} imply
\[
\Lambda=\Lambda_B+\imath\Lambda_C, \quad \text{or equivalently,} \quad \lambda_r=\beta_r+\imath\gamma_r, \quad 1\leq r\leq n.
\]
Thus, the eigenvalues of $B$ ($C$) are the real (imaginary) parts of the eigenvalues of $A$.

We will consider the following block Jacobi process on the normal matrix $A$,
\begin{equation}\label{jacobiABC}
A^{(k+1)}=\left(U^{(k)}\right)^*A^{(k)}U^{(k)}=\left(U^{(k)}\right)^*B^{(k)}U^{(k)}+\imath \left(U^{(k)}\right)^*C^{(k)}U^{(k)}, \quad k\geq0, \quad A^{(0)}=A,
\end{equation}
where $U^{(k)}$ are $\ubce_{\pi}(\varrho)$ matrices for some $0<\varrho \leq 1$. Here, for $k\geq 0$, $B^{(k)}$ ($\imath C^{(k)}$) is the Hermitian (skew-Hermitian) part of $A^{(k)}$. Obviously, $\lambda_r$, $\beta_r$, $\gamma_r$ are, for all $1\leq r\leq n$,
the eigenvalues of $A^{(k)}$, $B^{(k)}$, $C^{(k)}$, respectively. For each $k\geq 0$, the pivot submatrix $\widehat{U}_k$ will be chosen to diagonalize the pivot submatrix of $B^{(k)}$ (or $C^{(k)}$).

If $\pi=(1,1,\ldots,1)$ the block process (\ref{jacobiABC}) reduces to the element-wise Jacobi process
\begin{equation}\label{jacobiABC-el}
A^{(k+1)}=\left(R^{(k)}\right)^*A^{(k)}R^{(k)}=\left(R^{(k)}\right)^*B^{(k)}R^{(k)}+\imath \left(R^{(k)}\right)^*C^{(k)}R^{(k)}, \quad k\geq0, \quad A^{(0)}=A,
\end{equation}
where $R^{(k)}=R^{(k)}(i,j)$, $i=i(k)$, $j=j(k)$, are the complex plane rotations defined by the relations~\eqref{compl_rot_mat}--\eqref{simple_formulas}.
In this case, $R^{(k)}$ belongs to the class $\ubce_{\pi}(\varrho)$, for any $0<\varrho\leq1$, because the relation (\ref{angles}) implies $|\cos\phi_{ij}|=\cos\phi_{ij}\geq \sqrt{2}/2\geq\varrho(\sqrt{2}/2)$, for all $k\geq0$.

First, we consider the easier case when the real or imaginary parts of the eigenvalues of $A$ are simple.

\begin{Theorem}\label{tm:normal}
Let $A=B+\imath C\in\C^{n\times n}$ be a normal block matrix defined by the partition $\pi$. Let $(A^{(k)},k\geq0)$ be sequence of matrices generated by the iterative process~\eqref{jacobiABC}, where each transformation $U^{(k)}$ is chosen to diagonalize the pivot submatrix $\widehat{B}_{i(k)j(k)}^{(k)}$ of $B^{(k)}$. Let the cyclic pivot strategy be defined by $\calO\in\mathcal{B}_{sg}^{(m)}$. Let $0<\varrho \leq 1$ and let $U^{(k)}\in\ubce_{\pi}(\varrho)$, $k\geq 0$. Let the core algorithm be any convergent complex element-wise Jacobi method, and the core algorithm has attribute $\mathsf{R}$.  If $B$ has simple eigenvalues then $(A^{(k)},k\geq0)$ converges to a diagonal matrix of the eigenvalues of $A$, that is,
\[
\lim_{k\rightarrow\infty}A^{(k)}=\text{diag}(\lambda_1,\lambda_2,\ldots,\lambda_n).
\]
The same is true if each $U^{(k)}\in\ubce_{\pi}(\varrho)$ is chosen to diagonalize the pivot submatrix $\widehat{C}_{i(k)j(k)}^{(k)}$ of $C^{(k)}$, and if $C$ has simple eigenvalues.
\end{Theorem}

\begin{proof}
From (\ref{jacobiABC}) we obtain
$$B^{(k+1)}=\left(U^{(k)}\right)^*B^{(k)}U^{(k)}, \quad k\geq 0, \quad B^{(0)}=B,$$
which is a block Jacobi process on the Hermitian matrix $B$. Since all assumptions of the Theorem~\ref{tm:glconv} are satisfied, the assertions (i) and (ii) of that theorem imply
\begin{equation}\label{limB}
\lim_{k\rightarrow\infty}B^{(k)}=\Lambda_B=\diag(\beta_1,\beta_2,\ldots,\beta_n).
\end{equation}
Note that the matrices $B^{(k)}$ and $C^{(k)}$ commute for every $k\geq 0$. So, we have
\begin{equation}\label{BC=CB}
e_r^T(B^{(k)}C^{(k)})e_s=e_r^T(C^{(k)}B^{(k)})e_s, \quad 1\leq r,s\leq n, \quad k\geq 0,
\end{equation}
where $I_n=[e_1,\ldots ,e_n]$. We obtain
\begin{eqnarray}
  e_r^T(B^{(k)}C^{(k)})e_s = & (e_r^TB^{(k)})(C^{(k)}e_s)= \beta_r c_{rs}^{(k)}+\delta_{rs}^{(k)},
  \label{rBCs} \\
  e_r^T(C^{(k)}B^{(k)})e_s = & (e_r^TC^{(k)})(B^{(k)}e_s)= \beta_s c_{rs}^{(k)}+\zeta_{rs}^{(k)},
  \label{rCBs}
\end{eqnarray}
where
\begin{equation}\label{maxrs}
\max\{|\delta_{rs}^{(k)}|,|\zeta_{rs}^{(k)}|\} \leq \frac{\sqrt{2}}{2}\off (B^{(k)})\|C^{(k)}\|_2 = \frac{\sqrt{2}}{2}\off (B^{(k)})\|C\|_2\ \rightarrow0 \quad \text{as } k\rightarrow\infty.
\end{equation}

By subtracting the equations~\eqref{rBCs},~\eqref{rCBs}, and using~\eqref{BC=CB}, \eqref{maxrs} we obtain
\begin{equation}\label{limcrs}
\lim_{k\rightarrow\infty}c_{rs}^{(k)}(\beta_r-\beta_s)=0.
\end{equation}
Considering that all the eigenvalues of $B$ are distinct, $\beta_r\neq\beta_s$ implies
\[
\lim_{k\rightarrow\infty}c_{rs}^{(k)}=0 \quad \text{for}\ r\neq s.
\]
This proves $\lim_{k\rightarrow\infty}\off(C^{(k)})=0$. An easy calculation yields
\[
\off^2(A^{(k)}) = \off^2(B^{(k)}) +\off^2(C^{(k)}), \quad k\geq 0,
\]
and that proves
$$\lim_{k\rightarrow\infty}\off(A^{(k)})=0.$$

It remains to prove that in the later stage of the process the diagonal elements of $A^{(k)}$ cannot change their affiliation to the eigenvalues of $A$. This is easily solved by two arguments. First,
the eigenvalues $\lambda_r$, $1\leq r\leq n$, are also distinct. Namely, we have
\[
|\lambda_r-\lambda_s| \geq |\beta_r-\beta_s|, \quad 1\leq r\neq s\leq n.
\]
Second, since all $\beta_r$, $1\leq r\leq n$, are distinct, ultimately, all transformation matrices
$U^{(k)}$ approach the set of diagonal matrices. Therefore, for $k$ large enough, the change of $a_{ii}^{(k)}$ and $a_{jj}^{(k)}$ is so small that they cannot change their affiliation to the eigenvalues.

If each $U^{(k)}\in\ubce_{\pi}(\varrho)$ is chosen to diagonalize the pivot submatrix $\widehat{C}_{i(k)j(k)}^{(k)}$ of $C^{(k)}$, the proof is quite similar. Now, the roles of the matrices $B^{(k)}$ and $C^{(k)}$ in the proof are interchanged. Another way to prove the last assertion of the theorem is to consider the equation $-\imath A=C-\imath B$ and the iterative process $\imath (-A^{(k+1)})=\left(U^{(k)}\right)^*C^{(k)}U^{(k)} -\imath \left(U^{(k)}\right)^*B^{(k)}U^{(k)}$, $k\geq0$. The matrix $-\imath A$ and all matrices $-\imath A^{(k)}$, $k\geq 1$, are normal and $C$ and $C^{(k)}$, $k\geq 1$, are their Hermitian parts. It remains to apply the first assertion of the theorem to the
sequence $(-\imath A^{(k)},k\geq 0)$.
\end{proof}

It is recommended that, before the iterations start, the diagonal blocks of $B$ are diagonalized, such that in $B^{(0)}$, $B^{(0)}= U^{(0)*}BU^{(0)}$, the diagonal blocks are diagonal. Matrix $U^{(0)}$ is block diagonal, and we also have $A^{(0)}=\left(U^{(0)}\right)^*AU^{(0)}$.

For $\pi =(1,1,\ldots ,1)$ the following corollary holds.

\begin{Corollary}\label{cor:normal}
If $\pi =(1,1,\ldots,1)$, then Theorem~\ref{tm:normal} holds provided that each $U^{(k)}$ is the complex rotation.
\end{Corollary}

\begin{proof}
The proof follows the lines of the proof of Theorem~\ref{tm:normal}. Now, we consider the process~\eqref{jacobiABC-el} where each transformation matrix $R^{(k)}$ is defined by the relations~\eqref{compl_rot_mat}--\eqref{simple_formulas} for the matrix $B^{(k)}$ (or $C^{(k)}$ for the second assertion). The transformations are from the class $\ubce_{\pi}(\varrho)$ and the formulas~\eqref{simple_formulas} prevent the diagonal elements of $B^{(k)}$ from changing the affiliation in the later stage of the process. The assertion (iii) of Theorem~\ref{tm:normal} reduces to the assertion (ii). Since the eigenvalues of $B$ (or $C$) are simple, the complex rotations approach the set of diagonal matrices. Therefore, for $k$ large enough, the diagonal elements of $A^{(k)}$ cannot change affiliation.
\end{proof}

The case when neither $B$ nor $C$ has simple eigenvalues is much more complicated, especially when our goal is to construct a robust block Jacobi algorithm for solving the EVP of normal matrices. The solution would require a new research with a lot of numerical testing. Here, we give some hints.

Let us consider the case when each $U^{(k)}$ diagonalizes the pivot submatrix of $B^{(k)}$.
If $B$ has multiple eigenvalues then we would also require that ultimately the diagonal elements of $B^{(k)}$ are ordered non-increasingly (or non-decreasingly). Since we assume any generalized serial pivot strategy, this can be achieved by appropriately permuting the diagonal elements of $A^{(k)}$ after each cycle. In this way the limit matrix $\Lambda_B$ from the relation (\ref{limB}) will have
some equal successive diagonal elements, say
\[
\beta_1=\cdots =\beta_{s_1}>\beta_{s_1+1}=\cdots =\beta_{s_2}> \cdots > \beta_{s_{p-1}+1}=\cdots =\beta_{s_{p}},
\]
where $p$ is a number of different eigenvalues of $B$, and $\nu_1$, \ldots ,$\nu_p$,
\[
\nu_r = s_r-s_{r-1}, \quad 1\leq r\leq p, \ s_0=0, \ s_p=n,
\]
are multiplicities of the eigenvalues $\beta_{s_1}$,\ldots,$\beta_{s_p}$, respectively. We consider the case when at least one $\nu_r$ is larger than one. In such a case the sequence of Hermitian matrices $(C^{(k)},k\geq 0)$ may not converge to a diagonal Hermitian matrix.

Let $\nu = (\nu_1,\nu_2,\ldots ,\nu_p)$. We call $\nu$ \emph{natural partition} of $n$, and then we can call $\pi$ \emph{basic partition} of $n$. Using the block-matrix partition defined by $\nu$ we can write
\[
\Lambda_B = \left[\begin{array}{ccc}
\beta_1 I_{\nu_1} &  &     \\
\vdots & \ddots &\vdots \\
 & \cdots &  \beta_{s_p} I_{\nu_p}
        \end{array}\right],\ \
B^{(k)} = \left[\begin{array}{ccc}
\mathsf{B}_{11}^{(k)} & \cdots &  \mathsf{B}_{1p}^{(k)}   \\
\vdots & \ddots &\vdots \\
\mathsf{B}_{p1}^{(k)} & \cdots &  \mathsf{B}_{pp}^{(k)}
        \end{array}\right],\ \
C^{(k)} = \left[\begin{array}{ccc}
\mathsf{C}_{11}^{(k)} & \cdots &  \mathsf{C}_{1p}^{(k)}   \\
\vdots & \ddots &\vdots \\
\mathsf{C}_{p1}^{(k)} & \cdots &  \mathsf{C}_{pp}^{(k)}
        \end{array}\right]\begin{array}{c} \nu_1 \\ \vdots \\ \nu_p \end{array}.
\]
From the relation (\ref{limcrs}) we see that all off-diagonal blocks of the matrix $C^{(k)}$
tend to zero,
\[
\lim_{k\rightarrow\infty}\mathsf{C}_{rs}^{(k)} =0, \quad r\neq s.
\]
Suppose that after $K$ cycles of the block algorithm, we can detect gaps in the spectrum of $B$ and thus determine the partition $\nu$, and we can ensure that all off-diagonal blocks $\mathsf{B}_{rs}^{(KM)}$, $\mathsf{C}_{rs}^{(KM)}$, $r\neq s$ are $\epsilon$-small, where $\epsilon$ is much smaller than the minimum gap in the spectrum of $B$. The theory of nearly diagonal matrices (see \cite{Hari91}) implies that all $\off (\mathsf{B}_{rr}^{(KM)})$, $1\leq r\leq p$,  have to be $\epsilon^2$ small.

Then we diagonalize the diagonal blocks $\mathsf{C}_{rr}^{(KM)}$, $1\leq r\leq p$, using some convergent element-wise or block Jacobi method, depending on the size of the diagonal blocks. In this way, the quantities $\Off (B^{(k)})$ and $\Off (C^{(k)})$, where
\[
\Off (B^{(k)})=\|B^{(k)}-\mbox{diag}(\mathsf{B}_{11}^{(k)},\ldots ,\mathsf{B}_{pp}^{(k)})\|_F, \quad
\Off (C^{(k)})=\|C^{(k)}-\mbox{diag}(\mathsf{C}_{11}^{(k)},\ldots ,\mathsf{C}_{pp}^{(k)}\|_F,
\]
will not change for $k\geq KM$. Even more, since each $\beta_{s_r}I_{\nu_r}$ is invariant under unitary similarity transformations, all
$\off(\mathsf{B}_{rr}^{(k)})$, $1\leq r\leq p$, will remain to be $\epsilon^2$-small. This implies that at some cycle $K_1>K$,  all off-diagonal elements of $A^{(K_1M)}$ are $\epsilon$-small and those within the diagonal blocks are $\epsilon^2$-small.

In practical computation, if $\epsilon$ is chosen small enough, the stopping criterion could require one additional cycle. We can also pay attention whether the ultimate quadratic convergence of the process can occur and take it into consideration.

For theoretical considerations, one can chose positive numbers $\epsilon_1> \epsilon_2>\epsilon_3>\cdots $ such that $\lim_{t\rightarrow\infty}\epsilon_t=0$, and then, using the above analysis, find $K_1<K_2< \cdots $ such that $\off (A^{K_tM})\leq \epsilon_t$. This would ensure $\lim_{k\rightarrow\infty}\off (A^{(k)})=0$.

\section{The general Jacobi-type process for Hermitian matrices}\label{sec:Jacobi-type}

Let us consider a more general Jacobi-type block process for Hermitian matrices. The process uses the congruence transformation with nonsingular elementary block matrices.

Let $A\in\C^{n\times n}$ be a Hermitian matrix carrying the block-matrix partition defined by $\pi=(n_1,\ldots,n_m)$. We consider the iterative process of the form
\begin{equation}\label{jacobi-type}
A^{(k+1)}=\left(T^{(k)}\right)^*A^{(k)}T^{(k)}, \quad k\geq0, \quad A^{(0)}=A,
\end{equation}
where the transformations $T^{(k)}$ are nonsingular elementary block matrices
\[
T^{(k)}=\left[
        \begin{array}{ccccc}
          I &  &  &  &  \\
           & T_{ii}^{(k)} &  & T_{ij}^{(k)} &  \\
           &  & I &  &  \\
           & T_{ji}^{(k)} &  & T_{jj}^{(k)} &  \\
           &  &  &  & I \\
        \end{array}
      \right]\begin{array}{c}
                \\
               n_i \\
                \\
               n_j \\
                \\
             \end{array}.
\]
The process (\ref{jacobi-type}) has one nice property --- all  iteration matrices have to be Hermitian. We assume $A\neq0$, since, otherwise, all iterations have to be the null matrices.
Assumption $A\neq0$ implies that all $A^{(k)}\neq 0$, and $\|A^{(k)}\|>0$, $k\geq0$, holds for every matrix norm. Unlike in the case of the block Jacobi method, the matrices $T^{(k)}$ from~\eqref{jacobi-type} are not necessarily unitary. Moreover, it is not required that $T^{(k)}$ diagonalizes the pivot submatrix $\widehat{A}_{ij}^{(k)}$ of $A^{(k)}$.

For the process (\ref{jacobi-type}) we make the following assumptions.

\begin{itemize}
\item[\textbf{A1}] The pivot strategy is defined by any ordering $\calO\in\mathcal{B}_{sg}^{(m)}$.
\item[\textbf{A2}] There is a sequence of unitary elementary block matrices $(U^{(k)},k\geq0)$, such that
\[
\lim_{k\rightarrow\infty}(T^{(k)}-U^{(k)})=0.
\]
\item[\textbf{A3}] For the diagonal block $T_{ii}^{(k)}$ of $T^{(k)}$, we have
\[
\sigma =\liminf_{k\rightarrow\infty}\sigma^{(k)}>0, \quad \text{where} \ \sigma^{(k)}=\sigma_{\min}\big{(}T_{ii}^{(k)}\big{)}, \ k\geq 0.
\]
\end{itemize}

The assumption $\textbf{A3}$ is obviously satisfied if the matrices $U^{(k)}$ are from the class $\ubce_{\pi}(\varrho)$ for some $0<\varrho\leq1$. Since we assume that $\textbf{A2}$ holds, assumption $\textbf{A3}$ can also be written using the diagonal blocks $U_{ii}^{(k)}$, instead of $T_{ii}^{(k)}$.

The next theorem is a generalization of~\cite[Theorem 5.1]{Hari15}. The difference is in the assumption $\textbf{A1}$. In~\cite{Hari15} the pivot strategies are weakly equivalent to the column-cyclic strategy $I_{\mathcal{O}_{\text{col}}}$, while here we have $\calO\in\mathcal{B}_{sg}^{(m)}$. The proof is similar to the proof of~\cite[Theorem 5.1]{HaBe17} which considers real matrices.

\begin{Theorem}\label{tm:Jacobi-type}
Let $A\in\mathbb{C}^{n\times n}$, $A\ne0$, be a Hermitian block matrix defined by the partition $\pi=(n_1,\ldots,n_m)$, and let the sequence of matrices $(A^{(k)},k\geq0)$ be generated by the block Jacobi-type process~\eqref{jacobi-type}. If the assumptions $\mathbf{A1}-\mathbf{A3}$ are satisfied, then the following two conditions are equivalent:
\begin{itemize}
\item[(i)] $\displaystyle \lim_{k\rightarrow\infty} \frac{\emph{\off}(\widehat{A}_{i(k)j(k)}^{(k+1)})}{\|A^{(k)}\|_F}=0$,
\item[(ii)] $\displaystyle \lim_{k\rightarrow\infty} \frac{\emph{\off}(A^{(k)})}{\|A^{(k)}\|_F}=0$.
\end{itemize}
\end{Theorem}

To prove Theorem~\ref{tm:Jacobi-type} we need several lemmas.

\begin{Lemma}\label{lemma:Jacobi-type}
Let $(A^{(k)},k\geq0)$ be the sequence generated by the iterative process~\eqref{jacobi-type} and let $a^{(k)}=\ve(A^{(k)})$, $k\geq0$. Let the condition (i) of Theorem~\ref{tm:Jacobi-type} and the assumption $\mathbf{A2}$ hold. Then
\[
a^{(k+1)}=\mathcal{R}_{i(k)j(k)}a^{(k)}+g^{(k)}, \quad \mathcal{R}_{i(k)j(k)}\in\calRl_{ij}, \quad k\geq0,
\]
and
\[
\lim_{k\rightarrow\infty}\frac{g^{(k)}}{\|A^{(k)}\|_F}=0.
\]
\end{Lemma}

\begin{proof}
We can write the $k$th step of~\eqref{jacobi-type} as
\begin{equation}\label{lem61-0}
A^{(k+1)}=\left(U^{(k)}+(T^{(k)}-U^{(k)})\right)^*A^{(k)}(U^{(k)}+(T^{(k)}-U^{(k)}))=\left(U^{(k)}\right)^*A^{(k)}U^{(k)}+F^{(k)},
\end{equation}
where $U^{(k)}$ are unitary elementary block matrices matrices from the assumption $\textbf{A2}$.
We have
\begin{equation}\label{Fk}
F^{(k)}=\left(U^{(k)}\right)^*A^{(k)}(T^{(k)}-U^{(k)})+\left(T^{(k)}-U^{(k)}\right)^*A^{(k)}(T^{(k)}-U^{(k)})+\left(T^{(k)}-U^{(k)}\right)^*A^{(k)}U^{(k)}.
\end{equation}
We obtain
\[
\frac{\|F^{(k)}\|_F}{\|A^{(k)}\|_F} \leq 2\|T^{(k)}-U^{(k)}\|_F + \|T^{(k)}-U^{(k)}\|_F^2,
\]
and the assumption $\textbf{A2}$ implies
\begin{equation}\label{lem61-1}
\lim_{k\rightarrow\infty}\frac{\|F^{(k)}\|_F}{\|A^{(k)}\|_F}=0.
\end{equation}

Using the definition of the complex block Jacobi annihilator, for each $k$, we can write the mapping $A^{(k)}\mapsto \left(U^{(k)}\right)^*A^{(k)}U^{(k)}$ as
\begin{equation}\label{lem61-2}
a^{(k)}\mapsto \mathcal{R}_{i(k)j(k)}a^{(k)}+v^{(k)}, \quad \mathcal{R}_{i(k)j(k)}\in\calRl_{ij}, \quad k\geq0.
\end{equation}
Here, $v^{(k)}$ is a vector of the same size and the same partition as $a^{(k)}$, containing only  $2n_j$ non-zero blocks, those corresponding to the columns of $A_{i(k)j(k)}^{(k+1)}$ and rows of $A_{j(k)i(k)}^{(k+1)}$ (see the relation (\ref{vec})). Hence,
\[
\|v^{(k)}\|_2=\off(\widehat{A}_{i(k)j(k)}^{(k+1)}), \quad k\geq0,
\]
and, by the assumption of the lemma, that is the condition (i) of Theorem~\ref{tm:Jacobi-type}, we have
\begin{equation}\label{lem61-3}
\lim_{k\rightarrow\infty}\frac{\|v^{(k)}\|_2}{\|A^{(k)}\|_F}=0.
\end{equation}

Let
\begin{equation}\label{lem61-3a}
g^{(k)}=\ve(F^{(k)})+v^{(k)}, \quad k\geq 0.
\end{equation}
From the relations~\eqref{lem61-1} and~\eqref{lem61-3}, we obtain $\displaystyle \lim_{k\rightarrow\infty}\frac{\|g^{(k)}\|_2}{\|A^{(k)}\|_F}=0$. This, together with (\ref{lem61-0}) and (\ref{lem61-2}), proves the assertion of the lemma.
\end{proof}

\begin{Lemma}\label{lemma:Hari}
Let the assumptions of Lemma~\ref{lemma:Jacobi-type} hold. Then
\begin{itemize}
\item[(i)] $\displaystyle \lim_{k\rightarrow\infty}\frac{\|A^{(k+l)}\|_F}{\|A^{(k)}\|_F}=1, \quad l\geq0,$
\item[(ii)] $\displaystyle \lim_{k\rightarrow\infty}\frac{a^{(k)}}{\|A^{(k)}\|_F}=0 \quad \text{if and only if} \quad \lim_{t\rightarrow\infty}\frac{a^{(tM)}}{\|A^{(tM)}\|_F}=0$, \quad \text{where } $M=\frac{m(m-1)}{2}.$
\end{itemize}
\end{Lemma}

\begin{proof}
The proof is the same as that of~\cite[Lemma 5.3]{Hari15}.
\end{proof}

The next result is given in~\cite[Proposition 4.2]{BePe23}.

\begin{Lemma}[\cite{BePe23}]\label{lemma:sequence}
Let $(x_r,r\geq0)$ be a sequence of nonnegative real numbers such that
\[
x_{r+1} \leq \gamma x_r+c_r, \quad 0\leq\gamma<1.
\]
If \ \ $\lim_{r\to\infty}c_r=0$, then
\[
\lim_{r\to\infty}x_r=0.
\]
\end{Lemma}

Now, we are ready to prove the Theorem~\ref{tm:Jacobi-type}.

\begin{proof}[Proof of the Theorem~\ref{tm:Jacobi-type}]
Clearly, the condition $(ii)$ implies the condition $(i)$. We need to show that $(i)$ implies $(ii)$. Let us assume that the condition $(i)$ holds.

By the assumption \textbf{A1}, the pivot strategy is defined by  $\calO\in\mathcal{B}_{sg}^{(m)}$.
Let us consider the cycle $t$, $t\geq1$, of the iterative process~\eqref{jacobi-type}.
It consists of the steps $(t-1)M,(t-1)M+1,\ldots,tM-1$.
Let
$a^{[t]}=a^{(tM)}$, $A^{[t]}=A^{(tM)}$, $g^{[t]}=g^{(tM)}$, $t\geq1$, and $a^{[0]}=a^{(0)}$, $A^{[0]}=A^{(0)}$, $g^{[0]}=g^{(0)}$.

From the Lemma~\ref{lemma:Jacobi-type}, we conclude that during the cycle $t$ of the process~\eqref{jacobi-type} the off-diagonal elements satisfy the recursion
\begin{equation}\label{rel:Jtypecycle}
a^{[t]}=\calJ_t a^{[t-1]}+g^{[t]}, \quad \calJ_t\in\Jlo, \quad t\geq1,
\end{equation}
where
\[
g^{[t]}=g^{(tM-1)}+\sum_{\ell =0}^{M-2}\calR_{i(tM-1)j(tM-1)}\cdots
\calR_{i((t-1)M+\ell+1)j((t-1)M+\ell+1)}g^{((t-1)M+\ell)}, \quad t\geq 1.
\]
Since $\|\calR_{i(k)j(k)}\|_2=1$, $k\geq0$, we obtain
\begin{equation}\label{tm61-1}
\|g^{[t]}\|_2\leq \sum_{\ell =0}^{M-1}\|g^{((t-1)M+\ell)}\|_2, \quad t\geq 1.
\end{equation}
Using the inequality~\eqref{tm61-1} and Lemma~\ref{lemma:Hari}(i), it is easy to prove that the condition
$\lim_{k\rightarrow\infty}\frac{g^{(k)}}{\|A^{(k)}\|_F}=0$ implies
\begin{equation}\label{tm61-2}
\lim_{t\rightarrow\infty}\frac{g^{[t]}}{\|A^{[t]}\|_F}=0.
\end{equation}
Let
\[
b^{[t]}\coloneqq\frac{a^{[t]}}{\|A^{[t]}\|_F}, \quad t\geq1.
\]
By dividing the equation~\eqref{rel:Jtypecycle} with $\|A^{[t]}\|_F$ we get
\begin{align}
b^{[t]} & =\calJ_t\frac{a^{[t-1]}}{\|A^{[t]}\|_F} +\frac{g^{[t]}}{\|A^{[t]}\|_F} = \calJ_t b^{[t-1]}\frac{\|A^{[t-1]}\|_F}{\|A^{[t]}\|_F} +\frac{g^{[t]}}{\|A^{[t]}\|_F} \nonumber \\
& = \calJ_t b^{[t-1]} + c^{[t-1]},\quad t\geq1, \label{rel:Jtypeb}
\end{align}
where
\begin{equation}\label{tm61-3}
c^{[t-1]}=\calJ_t b^{[t-1]}\left(\frac{\|A^{[t-1]}\|_F}{\|A^{[t]}\|_F}-1\right)  +\frac{g^{[t]}}{\|A^{[t]}\|_F} \rightarrow 0, \quad \text{as } t\rightarrow\infty.
\end{equation}
Here we used the relation~\eqref{tm61-2}, Lemma~\ref{lemma:Hari}(i) and the inequality
$\|\calJ_t b^{[t-1]}\|_2\leq \| b^{[t-1]}\|_2 =1$, $t\geq1$.

Next, we use both assumptions \textbf{A2} and \textbf{A3}. Since $\sigma>0$, there is $t_0\geq1$
such that
\[
\| T^{(k)}-U^{(k)}\|_2 \leq \frac{1}{4}\sigma  \quad \text{and} \quad \sigma_{min}(T_{ii}^{(k)}) \geq \frac{3}{4}\sigma, \quad k\geq t_0M.
\]
From the perturbation theorem for the singular values we obtain, for any index $i$,
\begin{eqnarray*}
\sigma_{min}(U_{ii}^{(k)}) = &  \sigma_{min}(T_{ii}^{(k)}-(T_{ii}^{(k)}-U_{ii}^{(k)}))
  \geq \sigma_{min}(T_{ii}^{(k)})-\|T_{ii}^{(k)}-U_{ii}^{(k)}\|_2\\
\geq &  \frac{3}{4}\sigma - \|T^{(k)}-U^{(k)}\|_2 \geq  \frac{3}{4}\sigma - \frac{1}{4}\sigma
   = \frac{1}{2}\sigma, \quad k\geq t_oM.
\end{eqnarray*}
We set
\[
\varrho=\frac{1}{2}\sigma.
\]
Each $U^{(k)}$, $k\geq t_0M$, belongs to the class $\ubce_{\pi}(\varrho)$ of complex unitary elementary block matrices with partition $\pi$. Then, the block Jacobi operator $\calJ_{t}$, $t\geq t_0$, belongs to $\Jlo^{\ubce_{\pi}(\varrho)}$.

Without loss of generality, we can assume that the chain connecting the pivot orderings $\mathcal{O}\in\mathcal{B}_{sg}^{(m)}$ and $\mathcal{O''}\in\mathcal{B}_{sp}^{(m)}$ is in the canonical form and contains $d$ shift equivalences. We observe any $d+1$ successive cycles of the process~\eqref{rel:Jtypeb} after the cycle $t_0$ has been completed. We have
\begin{equation}\label{tm61-4}
b^{[t+d+1]}=\calJ_{t+d+1}\cdots\calJ_{t+2}\calJ_{t+1} b^{[t]}+h^{[t]}, \quad t\geq t_0,
\end{equation}
where
$$h^{[t]}=\calJ_{t+d+1}\cdots\calJ_{t+2}c^{[t]} + \calJ_{t+d+1}\cdots\calJ_{t+3}c^{[t+1]} + \cdots + \calJ_{t+d+1}c^{[t+d-1]} + c^{[t+d]}.$$
Note that $\|\calJ_t\|_2\leq1$, $t\geq t_0$. Since the relation (\ref{tm61-3}) holds, using a similar analysis as the one that resulted in relation (\ref{tm61-2}), we obtain
$$\lim_{t\rightarrow\infty}h^{[t]}=0.$$
Applying the vector norm to (\ref{tm61-4}), and using the Theorem~\ref{tm:sg_jop} to bound $\|\calJ_{t+d+1}\cdots\calJ_{t+1}\|_2$, we obtain
\begin{equation}\label{rel:sequenceb}
\|b^{[t+d+1]}\|_2\leq \mu_{\pi,\varrho} \|b^{[t]}\|_2+\|h^{[t]}\|_2, \quad t\geq t_0,
\end{equation}
where $0\leq\mu_{\pi,\varrho}<1$. Let $y_r=\|b^{[r(d+1)]}\|_2$, $z_r=\|h^{[r(d+1)]}\|_2$, $r\geq r_0$, where $r_0$ is such that $r_0(d+1)\geq t_0$. Then, the relation~\eqref{rel:sequenceb} implies that
\begin{equation}\label{tm61-6}
y_{r+1}\leq \mu_{\pi,\varrho} y_r+z_r, \quad r\geq r_0,
\end{equation}
holds, with $\lim_{r\rightarrow\infty}z_r=0$. Applying Lemma~\ref{lemma:sequence} to the sequence
$(y_r,r\geq r_0 )$ generated by the rule (\ref{tm61-6}) one obtains $\lim_{r\rightarrow\infty}y_r=0$.
Thus, we have
\begin{equation}\label{tm61-7}
\lim_{r\rightarrow\infty}b^{[r(d+1)]}=0.
\end{equation}
To show that $\lim_{t\rightarrow\infty}b^{[t]}=0$, let $r(d+1)\leq t< (r+1)(d+1)$, $r\geq r_0$. We can write
$t=r(d+1)+\xi$ for some integer $\xi$, $0\leq\xi <d+1$. From the relation~\eqref{rel:Jtypeb} we get
\[
\|b^{[t]}\|_2 = \|\calJ_{t}\cdots \calJ_{r(d+1)+1} b^{[r(d+1)]}+\tilde{c}_{t,\xi}\|_2
\leq \|b^{[r(d+1)]}\|_2+\|\tilde{c}_{t,\xi}\|_2 ,
\]
where
\begin{equation}\label{tm61-8}
\|\tilde{c}_{t,\xi}\|_2 \leq \|c^{[t-1]}\|_2+ \|c^{[t-2]}\|_2+\cdots\ + \|c^{[t-\xi]}\|_2.
\end{equation}
From the relations~\eqref{tm61-8} and~\eqref{tm61-3} we obtain $\lim_{r\rightarrow\infty}\|\tilde{c}_{t,\xi}\|_2 =0$. Together with~\eqref{tm61-7}, this implies
\begin{equation}\label{tm61-9}
\lim_{t\rightarrow\infty}b^{[t]}=0.
\end{equation}

Now, we can use Lemma~\ref{lemma:Hari}(ii) to obtain $\lim_{k\rightarrow\infty}b^{(k)}=0$, and we are done. However, in order to give the reader a better understanding, we provide a brief proof for this claim.

The sequence $(b^{[t]},t\geq t_0)$ that converges to zero is a subsequence of the sequence $(b^{(k)},k\geq t_0M)$.
Let $k$ be such that $tM\leq k<(t+1)M$, $t\geq t_0$. For some integer $\eta_k$, $0\leq\eta_k<M$, we have $k=(t-1)M+\eta_k$. From Lemma~\ref{lemma:Jacobi-type} we get
\begin{align*}
a^{(k)} & = \mathcal{R}_{i(k-1)j(k-1)}a^{(k-1)}+g^{(k-1)}
=\mathcal{R}_{i(k-1)j(k-1)}\left(\mathcal{R}_{i(k-2)j(k-2)}a^{(k-2)}+g^{(k-2)}\right)+g^{(k-1)}\\
& = \cdots = \mathcal{R}_{i(k-1)j(k-1)}\cdots \mathcal{R}_{i(tM)j(tM)}a^{(tM)}+g^{(k-1)}
   +\mathcal{R}_{i(k-1)j(k-1)}g^{(k-2)}+\cdots \\
& \quad + \mathcal{R}_{i(k-1)j(k-1)}\cdots \mathcal{R}_{i(tM+1)j(tM+1)}g^{(tM)}.
\end{align*}
Applying the norm to the left and right side we obtain
\[
\|a^{(k)}\|_2\leq  \|a^{(tM)}\|_2+ \|g^{(k-1)}\|_2+\cdots +\|g^{(tM)}\|_2, \quad k=(t-1)M+\eta_k.
\]
Here we used the fact that the spectral norm of the Jacobi annihilator equals one.
Using the Lemma~\ref{lemma:Hari}, and the relations~\eqref{tm61-9} and~\eqref{tm61-3}, we have
\begin{eqnarray*}
\|b^{(k)}\|_2=\frac{\|a^{(k)}\|_2}{\|A^{(k)}\|_F}&\leq& \|b^{(tM)}\|_2\ \frac{\|A^{(tM)}\|_F}{\|A^{(tM+1)}\|_F}\cdots  \frac{\|A^{(k-1)}\|_F}{\|A^{(k)}\|_F} +\frac {\|g^{(k-1)}\|_2}{\|A^{(k-1)}\|_F}\ \frac{\|A^{(k-1)}\|_F}{\|A^{(k)}\|_F}+\cdots \\
&& +\frac{\|g^{(tM)}\|_2}{\|A^{(tM)}\|_F}\ \frac{\|A^{(tM)}\|_F}{\|A^{(tM+1)}\|_F}\cdots  \frac{\|A^{(k-1)}\|_F}{\|A^{(k)}\|_F} \rightarrow 0, \quad \text{as } k\rightarrow\infty.
\end{eqnarray*}
This proves the theorem since $\off(A^{(k)})=\|a^{(k)}\|_2$, for all $k$.
\end{proof}

Similarly as in~\cite{HaBe17}, we have two corollaries that are very useful in applications.

\begin{Corollary}\label{cor:5.3}
Theorem~\ref{tm:Jacobi-type} holds provided the relation~\eqref{jacobi-type} is replaced by
\begin{equation}\label{jct1}
A^{(k+1)}=\left(T^{(k)}\right)^*A^{(k)}T^{(k)}+E^{(k)}, \quad k\geq0, \quad A^{(0)}=A,
\end{equation}
where
\begin{equation}\label{jct1a}
\lim_{k\rightarrow\infty}\frac{\emph{\off}(E^{(k)})}{\|A^{(k)}\|_F}=0 \quad \text{and} \quad
E^{(k)}\neq \left(T^{(k)}\right)^*A^{(k)}T^{(k)}, \quad k\geq 0.
\end{equation}
The condition (\ref{jct1a}) can be replaced by the requirement that $E^{(k)}=0$ whenever
$A^{(k)}=0$, for some $k$.
\end{Corollary}

\begin{proof}
By inspecting the proof of the Theorem~\ref{tm:Jacobi-type} we see that it suffices to update the definition of the matrix $F^{(k)}$ from the relation~\eqref{Fk} (by adding the matrix $E^{(k)}$) and the vector $g^{(k)}$ from the proof of Lemma~\ref{lemma:Jacobi-type}. In the relation~\eqref{lem61-3a} the vector $g^{(k)}$ should now include $\ve(E^{(k)})$, i.e.,
\begin{equation}\label{lem61-3b}
g^{(k)}=\ve(F^{(k)})+v^{(k)}+\ve(E^{(k)}), \quad k\geq 0.
\end{equation}
\end{proof}

\begin{Corollary}\label{cor:5.4}
Let $A\in\mathbb{C}^{n\times n}$, $A\ne0$, be a Hermitian block matrix defined by the partition $\pi=(n_1,\ldots,n_m)$, and let the sequence of matrices $(A^{(k)},k\geq0)$ be generated by the block Jacobi-type process~\eqref{jct1}. Let the assumptions $\mathbf{A1}-\mathbf{A3}$ hold. If the sequence $(A^{(k)},k\geq 0)$ is bounded and
$$\lim_{k\rightarrow\infty}\emph{\off}(E^{(k)})=0,$$
then the following two conditions are equivalent:
\begin{itemize}
\item[(iii)] $\displaystyle \lim_{k\rightarrow\infty} \emph{\off}(\widehat{A}_{i(k)j(k)}^{(k+1)})=0$,
\item[(iv)] $\displaystyle \lim_{k\rightarrow\infty} \emph{\off}(A^{(k)})=0$.
\end{itemize}
\end{Corollary}

\begin{proof}
The proof is quite similar to that of~\cite[Corollary~5.3]{HaBe17}, so we provide its brief version.
Instead of the relations~\eqref{lem61-0} and~\eqref{lem61-3a} we have $A^{(k+1)}=\left(U^{(k)}\right)^*A^{(k)}U^{(k)}+F^{(k)}+E^{(k)}$ and~\eqref{lem61-3b}. The assumption \textbf{A2} and the assumption on the sequence $(A^{(k)},k\geq 0)$ yield
\[
\|F^{(k)}\|_2 \leq \sup \left\{\|A^{(k)}\|_2|,k\geq0 \right\}\left(2\|T^{(k)}-U^{(k)}\|_F+\|T^{(k)}-U^{(k)}\|_F^2\right) \rightarrow 0, \quad \text{as } k\rightarrow\infty .
\]
Following the proofs of Lemma~\ref{lemma:Jacobi-type} and Theorem~\ref{tm:Jacobi-type}, we easily conclude that $\lim_{k\rightarrow\infty} g^{(k)}=0$, $\lim_{t\rightarrow\infty} g^{[t]}=0$,
$\lim_{t\rightarrow\infty} a^{[t]}=0$ and finally $\lim_{k\rightarrow\infty} a^{(k)}=0$.
\end{proof}

We end the theoretical considerations with one important application of Corollary~\ref{cor:5.4} in the next section.
Another  applications of the presented theory lie in the area of complex block Jacobi methods for the positive definite generalized eigenvalue problem $Ax=\lambda Bx$, with complex Hermitian $A$ and complex positive definite $B$. Here, one can use the complex HZ method from \cite{Hari19,Hari22} or the complex CJ method from \cite{Hari18} as the core algorithm. The convergence proof will use the theory from \cite{Hari21} together with Theorem~\ref{tm:Jacobi-type} and its corollaries.

\section{The global convergence of the complex block $J$-Jacobi method}\label{sec:J-Jacobi}

Here we briefly show how the Corollary~\ref{cor:5.4} can be used to prove the global convergence of the block $J$-Jacobi method for solving the $J$-Hermitian eigenvalue problem $Ax=\lambda Jx$, $x\neq0$. Matrix $A$ is complex Hermitian positive definite and $J=\diag(I_{\nu} ,-I_{n-\nu})$, for some $1\leq \nu <n$. This eigenvalue problem arises when the goal is to compute the eigenvalues and eigenvectors of an indefinite Hermitian matrix to high relative accuracy (see~\cite{Slap92,Veselic00,Slap03,SlapTruh03}).

To solve the $J$-Hermitian eigenvalue problem one can use the complex block $J$-Jacobi method which is described and analyzed in~\cite{HaSiSi14}. The complex block $J$-Jacobi method is a generalization of the known element-wise $J$-Jacobi method which was derived by Veseli\'{c} in~\cite{Veselic93}. As it is proven in~\cite{HaSiSi14} the method computes the eigenvalues of the matrix pair $(A,J)$ to high relative accuracy when the standard conditions are met. In the same paper the global convergence of the method is proven under the cyclic pivot strategies that are weakly equivalent to the row-cyclic strategy. Here we show that the global convergence result can be extended to hold for any generalized serial strategy.

The block method is defined by the partition $\pi =(n_1,\ldots ,n_m)$ which has to be subpartition of the partition $(\nu ,n-\nu)$. We can write $\nu = n_1+\cdots +n_p$ and $n-\nu =n_{p+1}+\cdots +n_{p+q}$. It implies that $m=p+q$. The method uses $J$-unitary elementary block matrices $T^{(k)}$ which satisfy $\left(T^{(k)}\right)^*JT^{(k)}= J$, $k\geq 0$. The iterative process has the form
$$A^{(k+1)}=\left(T^{(k)}\right)^*A^{(k)}T^{(k)}, \quad k\geq 0, \quad A^{(0)}=A.$$
As before, it is advantageous for the algorithm if the diagonal blocks of the starting matrix $A^{(0)}$ are diagonal. This can be achieved using an unitary block diagonal matrix $U^{(0)}$ in the same way as it is done for the block Jacobi method for Hermitian matrices. Thus, we can assume $A^{(0)}=U^{(0)*}AU^{(0)}$.

The transformation matrices $T^{(k)}$ are either unitary elementary block matrices (provided that $1\leq i(k)<j(k)\leq p$ or $p+1\leq i(k)<j(k)\leq m$) or they are hyperbolic elementary block matrices (provided that $1\leq i(k)\leq p< j(k)\leq m$). If $T^{(k)}$ is unitary, we will assume that it belongs to the class $\ubce_{\pi}(\varrho)$ with $\varrho=1$.
Such transformation matrices will be called $\ubce_{\pi}$ unitary elementary block matrices.

\subsection{The global convergence of the element-wise $J$-Jacobi method}

Since the element-wise $J$-Jacobi method with $\pi =(1,1,\ldots ,1)$ is used as the core algorithm of the block method, let us recall the main formulas (see \cite[relation (2.10)]{HaSiSi10}). The transformations are complex unitary or hyperbolic rotations. We will display their pivot submatrices.  Using $i$, $j$ instead of $i(k)$, $j(k)$, respectively, and omitting $k$, we have
\[
\widehat{T} = \widehat{\Phi}\widehat{R}, \quad \widehat{\Phi}=\diag(e^{\imath \phi},1), \quad \phi=\arg (a_{ij}),
\]
and
\[
\widehat{R} = \left\{
\begin{array}{l}
 \left[\begin{array}{cc}
     \cosh\theta & \sinh\theta \\
     \sinh\theta & \cosh\theta
 \end{array}\right],
 \quad \tanh\theta=\frac{-2|a_{ij}|}{a_{ii}+a_{jj}}
 \text{ if } 1\leq i\leq \nu <j\leq n\\
\left[\begin{array}{cc}
   \  \cos\theta & \ \ \sin\theta \\
   \  \sin\theta & \ \ \cos\theta
\end{array}\right], \quad \tan\theta=\frac{2|a_{ij}|}{a_{ii}-a_{jj}} \text{ elsewhere.}
\end{array}
     \right.
\]
Here, the first component of the transformation $\widehat{\Phi}$ transforms the pivot submatrix $\widehat{A}$ into $|\widehat{A}|$. It transforms $a_{ij}$ and
$a_{ji}$ into $|a_{ij}|$ and $|a_{ji}|$. The transformation formulas for the diagonal elements are given either by the relation (\ref{simple_formulas}) (with angle $\phi$ replaced by $\theta$)
or by
\[
a_{ii}' = a_{ii}+\tanh\theta|a_{ij}|, \quad a_{ii}' = a_{jj}+\tanh\theta|a_{ij}|.
\]
The former (letter) formulas are used when the transformation is trigonometric (hyperbolic) complex rotation.

\begin{Theorem}\label{tm:glc-J-Jacobi-el}
Suppose the element-wise complex $J$-Jacobi method uses the cyclic pivot strategy defined by $\mathcal{O}\in\mathcal{B}_{sg}^{(m)}$. Then it is globally convergent.
\end{Theorem}

\begin{proof}
It is proven in~\cite[relation (2.14)]{HaSiSi10} that $\|A^{(k)}\|_F\leq \mbox{trace}(A)$, $k\geq0$, holds for the element-wise $J$-Jacobi method. It implies that the sequence $(A^{(k)},k\geq 0)$ is bounded.
Let us show that the assumptions \textbf{A1}--\textbf{A3} hold. The first one \textbf{A1} is just the selection of the pivot strategy which is the assumption of the theorem. The second one
\textbf{A2} is automatically fulfilled for the complex rotations. In the hyperbolic steps we have
(see~\cite{HaSiSi10})
$a_{i(k)j(k)}^{(k)}\rightarrow 0$ and $\theta_k\rightarrow 0$, as $k\rightarrow\infty$. This shows that $\|T^{(k)}-\widehat{\Phi}^{(k)}\|_2\rightarrow 0$, as $k\rightarrow\infty$, over the set of hyperbolic
steps. Here, $\widehat{\Phi}_k$ is unitary. In order to check \textbf{A3}, note that $\cos\theta_k\geq \sqrt{2}/2$, for unitary steps, and $\cosh\theta_k\geq 1$, for hyperbolic steps. In addition,
$|e^{\imath \phi_k}|=1$, for all $k$. Hence, $\sigma_{min}(T^{(k)}_{ii})\geq \sqrt{2}/2$, for all $k$.
Thus, all three assumption \textbf{A1}--\textbf{A3} hold for the element-wise $J$-Jacobi method. Since the pivot element is annihilated at each step, the condition (i) also holds. Using the Corollary~\ref{cor:5.4} with $E^{(k)}=0$, we obtain $\displaystyle\lim_{k\rightarrow\infty} \off(A^{(k)})=0$, for all $k$. The proof that $A^{(k)}$ converges to the diagonal matrix $\Lambda$ of the eigenvalues of $(A,J)$ is the same as in the proof of \cite[Theorem~3.7]{HaSiSi10}.
\end{proof}

When using the element-wise $J$-Jacobi method as the core algorithm for the block method, we actually invoke the Theorem~\ref{tm:glc-J-Jacobi-el} in each hyperbolic step. In the $k$th step of the block method the partition $(n_{i(k)},n_{j(k)})$ is used instead of $(\nu ,n-\nu)$.

\subsection{The global convergence of the block $J$-Jacobi method}

The block $J$-Jacobi method was researched in~\cite{HaSiSi14}. This facilitates our convergence proof because we can use here several results from \cite{HaSiSi14}. Recall that the unitary transformation is used when $1\leq i(k)<j(k)\leq p$ or $p+1\leq i(k)<j(k)\leq m$.

\begin{Theorem}\label{tm:J-Jacobi}
Let $A\in\mathbb{C}^{n\times n}$ be a Hermitian positive definite block matrix defined by the partition $\pi=(n_1,\ldots,n_m)$ and let $J=\mbox{diag}(I_{\nu},-I_{n-\nu})$. Let the $J$-Hermitian block Jacobi method be applied to the pair $(A,J)$ generating sequence of matrices $(A^{(k)},k\geq 0)$.  Let the unitary transformations belong to the class $\ubce_{\pi}$. If the pivot strategy is generalized serial one, then we have
\begin{equation*}\label{glconv_blockJ}
    \lim_{k\rightarrow\infty}\emph{\off} (A^{(k)})=0.
\end{equation*}
\end{Theorem}

\begin{proof}
In \cite[Section~3]{HaSiSi14} it is shown that $\|A^{(k)}\|_F\leq \mbox{trace}(A^{(k)})\leq \mbox{trace}(A)$, $k\geq0$. This proves that the sequence $(A^{(k)},k\geq0)$ is bounded.
Therefore, we can apply Corollary~\ref{cor:5.4} with the matrices $E^{(k)}$, $k\geq 0$, set to zero.
The $J$-Jacobi method diagonalizes the pivot submatrix at each step. Thus, the condition (iii) of
Corollary~\ref{cor:5.4} is fulfilled. The rest of the proof reduces to checking validity of the assumptions \textbf{A1} -- \textbf{A3}.

The first assumption \textbf{A1} is presumed.
The rest of the proof is the same as the appropriate part of the proof of \cite[Proposition~3.2]{HaSiSi14}. It uses the theory developed in \cite[Section~3]{HaSiSi14} for the block hyperbolic steps, and holds for any pivot strategy.
 \end{proof}

The following corollary completes the global convergence consideration for the block $J$-Jacobi method.

Recall the property $\mathsf{R}$ of the core algorithm in the block Jacobi method for Hermitian matrices. We will also use it in the block $J$-Jacobi method. If the transformation matrix is not unitary, we can still apply the QR factorization with column pivoting to $[T_{ii}T_{ij}]$, and obtain the permutation matrix in the form $\widehat{P}=I_{11'}\cdots I_{n_in_i'}$. Then, one can apply the procedure from  Definition~\ref{attr_R} to obtain $\widetilde{\widehat{P}}$. In this way the core algorithm in the block $J$-Jacobi method can get attribute $\mathsf{R}$.

\begin{Corollary}\label{cor:J-Jacobi}
Let all assumptions of Theorem~\ref{tm:J-Jacobi} hold.  Then the assertions (i)--(ii) hold.
\begin{itemize}
\item[(i)] If the eigenvalues of the pair $(A,J)$ are simple and the core algorithm has attribute $\mathsf{R}$, then $\Lambda_A=\lim_{k\rightarrow\infty} A^{(k)}$ exists, and $J\Lambda_A$ is the diagonal matrix of the eigenvalues of the pair $(A,J)$.
\item[(ii)] If the core algorithm has attribute $\mathsf{R}$ and if
  there is $k_0\geq 0$ such that for $k\geq k_0$ all diagonal elements affiliated with a multiple eigenvalue are contained within the same diagonal block, then $\Lambda_A=\lim_{k\rightarrow\infty} A^{(k)}$ exists, and $J\Lambda_A$ is the diagonal matrix of the eigenvalues of $(A,J)$.
\end{itemize}
\end{Corollary}

\begin{proof}
By the Theorem~\ref{tm:J-Jacobi} the block method converges to diagonal form.

\begin{itemize}
\item[(i)] The proof is very similar to the proof of Theorem~\ref{tm:glconv}(ii). One can also use the following observation. For $k$ large enough, $T^{(k)}$ is arbitrary close to some diagonal unitary matrix $\Psi^{(k)}$. Writing $T^{(k)}=\Psi^{(k)}+\Xi^{(k)}$, where $\Xi^{(k)}\rightarrow 0$ as $k\rightarrow\infty$, we easily conclude that, for $k$ large enough, the change of each diagonal element of $A^{(k)}$ is arbitrary small. It means that $JA^{(k)}$ converges to the diagonal matrix of the eigenvalues of $(A,J)$.
\item[(ii)] For the unitary transformations the proof follows the lines of the proof of Theorem~\ref{tm:glconv}(iii). Only at one place, the appearance of $\widehat{U}_{ij}$ should be replaced by $\widehat{T}_{ij}$. Now, consider the hyperbolic steps.
As $k$ increases over the set of hyperbolic steps, the transformation matrices $T^{(k)}$ approach the set of diagonal unitary matrices. Therefore, the hyperbolic transformations will not change affiliation of the diagonal elements in the later stage of the process.
\end{itemize}
\end{proof}

\section{Numerical experiments}\label{sec:num}

In the final section we give the results of our numerical experiments. We test the convergence speed and the accuracy of the block Jacobi method on different block matrices, with different block sizes and under different pivot strategies $I_{\calO}$, $\calO\in\mathcal{B}_{sp}^{(m)}$. Numerical tests are performed in MATLAB R2023b.

\begin{figure}[!h]
\begin{subfigure}{0.45\textwidth}
    \centering
    \includegraphics[width=\textwidth]{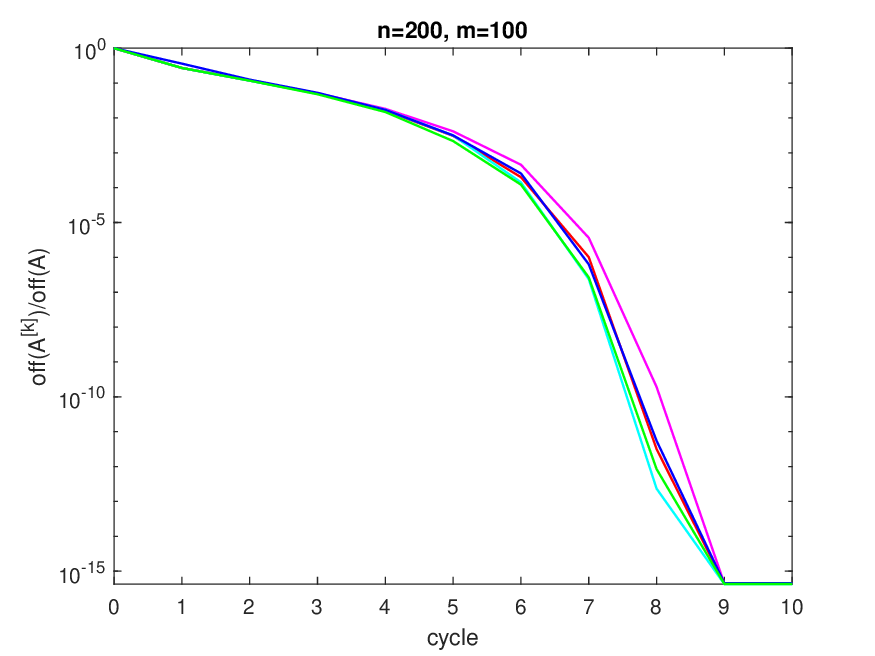}
    \caption{$2\times2$ blocks}
\end{subfigure}
\hfill
\begin{subfigure}{0.45\textwidth}
        \includegraphics[width=\textwidth]{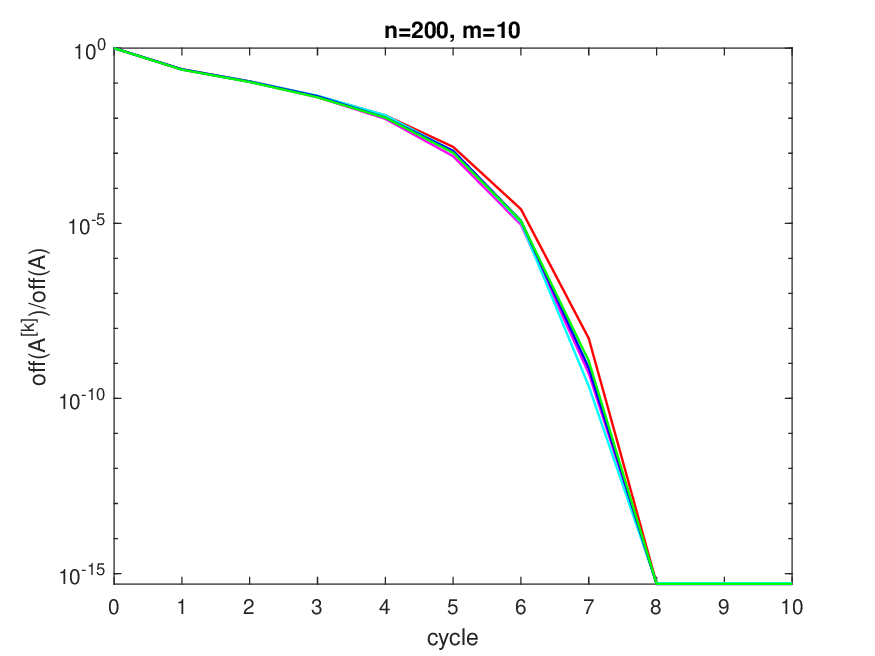}
    \caption{$20\times20$ blocks}
\end{subfigure}
\caption{Convergence of the off-norm on a random Hermitian block matrix $A\in\C^{200\times200}$ under five different generalized serial pivot strategies.}
\label{fig:strategies}
\end{figure}

In the first set of experiments we compare five pivot strategies $I_{\calO}$, $\calO\in\mathcal{B}_{sp}^{(m)}$. The block strategies are derived from the row-cyclic pivot strategy using random permutations in each row. The core algorithm is the elemet-wise Jacobi method under a generalized serial pivot strategy, randomly chosen and different in each block iteration. We work on a random Hermitian matrix $A\in\C^{200\times200}$ and use different block sizes. For $2\times2$ blocks, we have $m=100$ and the partition $\pi$ from~\eqref{partition} is such that $n_1=n_2=\cdots=n_{100}=2$. For $20\times20$ blocks, $m=10$ and $n_1=n_2=\cdots=n_{10}=20$. Results are presented in Figure~\ref{fig:strategies}. It can be observed that all strategies behave approximately in the same way.

\begin{figure}[!h]
\begin{subfigure}{0.45\textwidth}
    \centering
    \includegraphics[width=\textwidth]{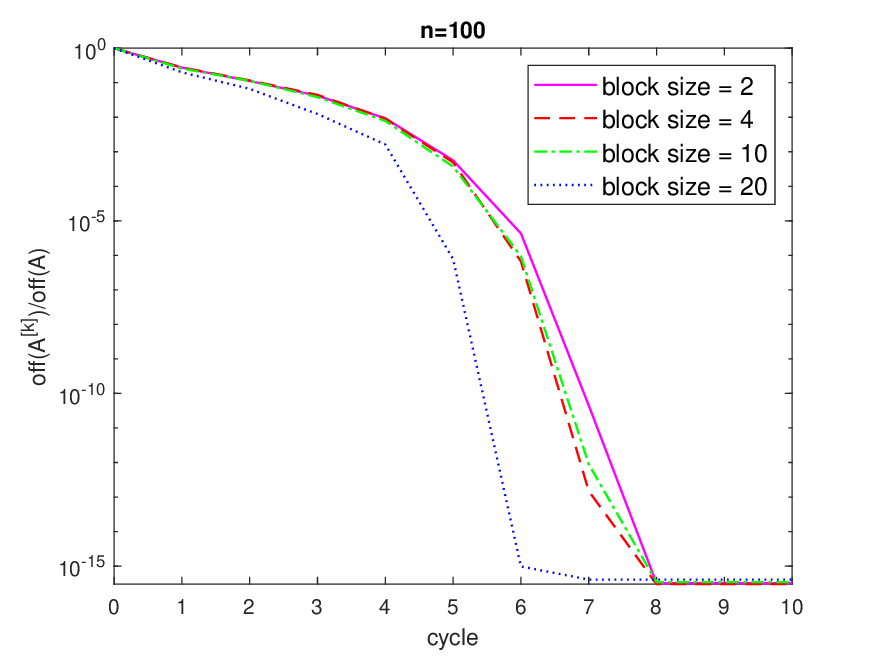}
    \caption{Random $100\times100$ Hermitian matrix}
\end{subfigure}
\hfill
\begin{subfigure}{0.45\textwidth}
        \includegraphics[width=\textwidth]{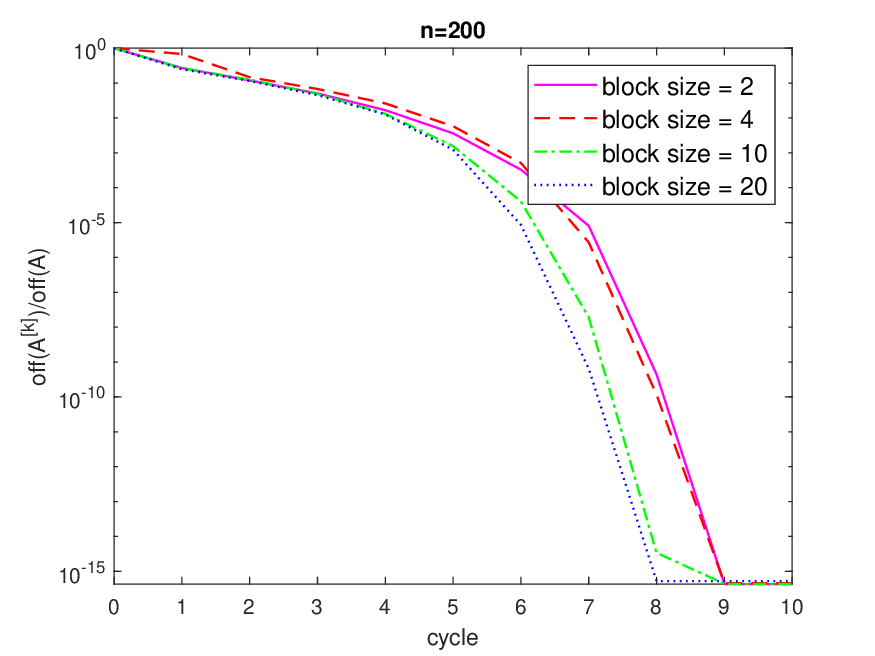}
    \caption{Random $200\times200$ Hermitian matrix}
\end{subfigure}
\caption{Convergence of the off-norm on a random Hermitian block matrix $A\in\C^{n\times n}$ with different block sizes.}
\label{fig:blocks}
\end{figure}

In Figure~\ref{fig:blocks} we compare the convergence of the off-norm for different block sizes on a random Hermitain matrix $A\in\C^{n\times n}$. The pivot strategy $I_{\calO}$, $\calO\in\mathcal{B}_{sp}^{(m)}$, is fixed and randomly chosen in the same way as in the previous example. Here, one can see that, in general, for the larger blocks, the algorithm needs less cycles to converge.

\begin{figure}[!h]
\begin{subfigure}{\textwidth}
    \centering
    \includegraphics[width=\textwidth]{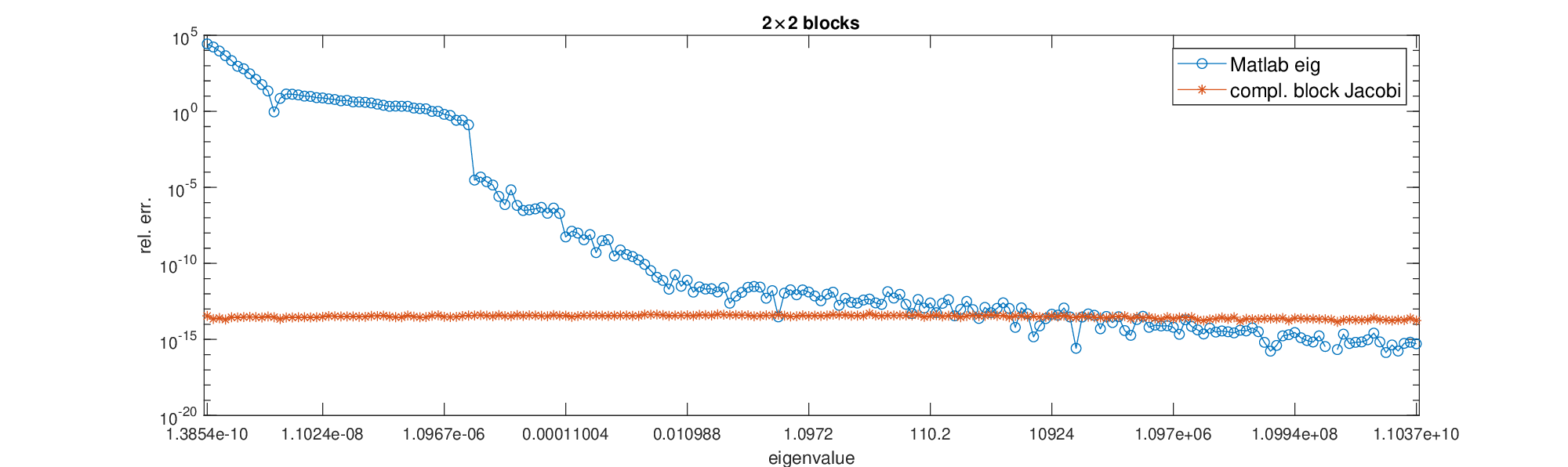}
\end{subfigure}
\begin{subfigure}{\textwidth}
    \centering
    \includegraphics[width=\textwidth]{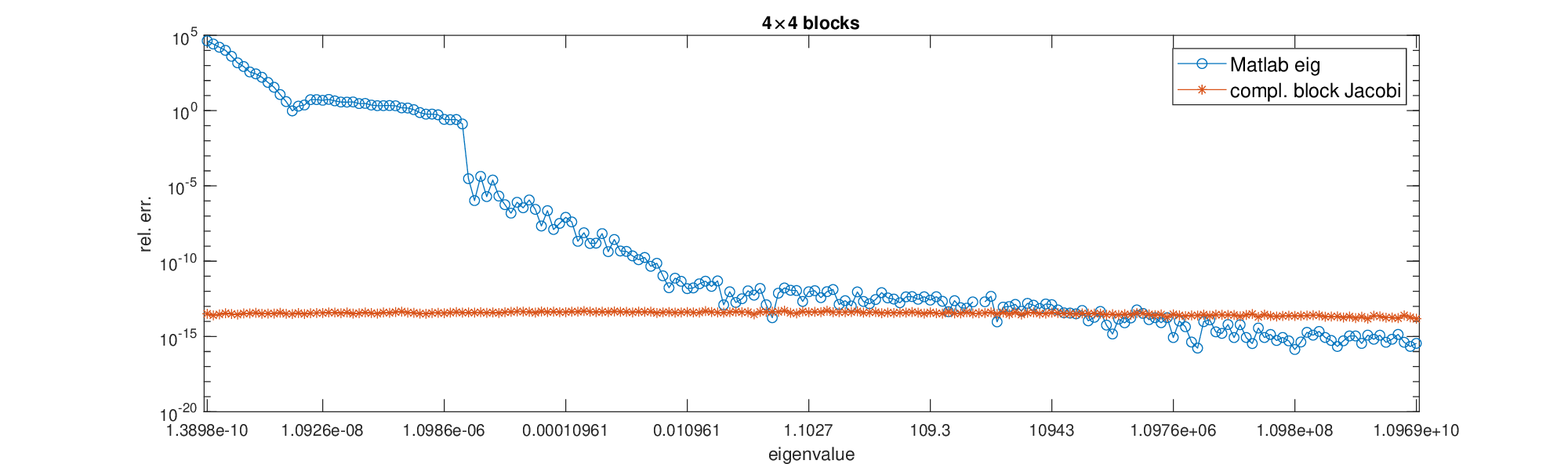}
\end{subfigure}
\begin{subfigure}{\textwidth}
    \centering
    \includegraphics[width=\textwidth]{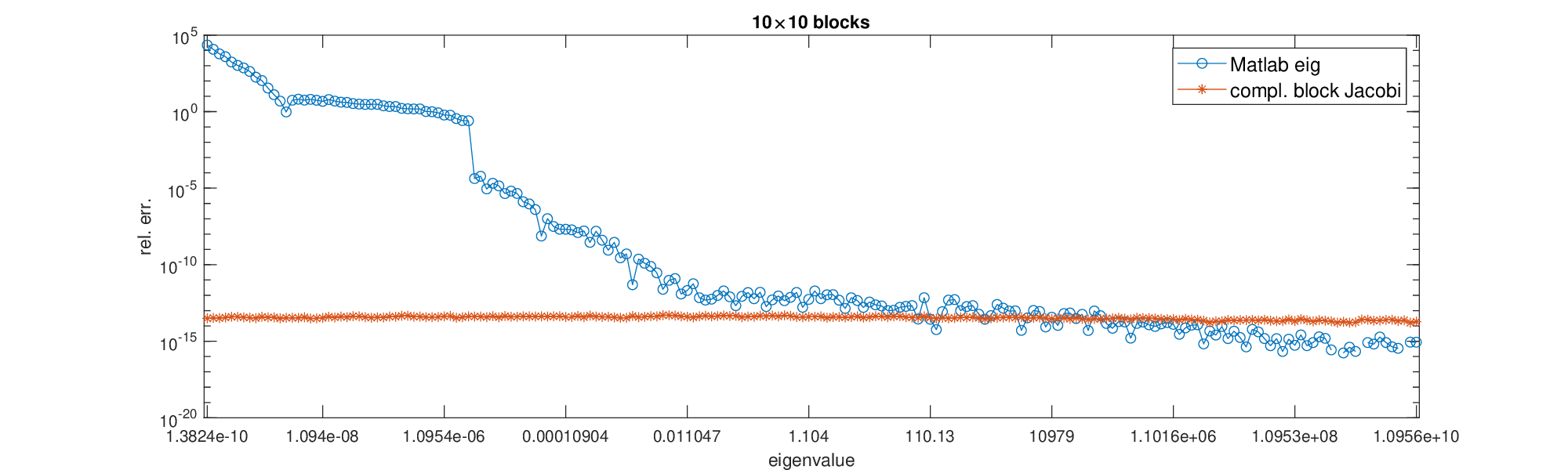}
\end{subfigure}
\caption{Accuracy of the computed eigenvalues on an ill-conditioned Hermitian block matrix $A\in\C^{200\times200}$ with different block sizes.}
\label{fig:accuracy}
\end{figure}

In Figure~\ref{fig:accuracy} we test the accuracy of the complex block Jacobi algorithm compared to the MATLAB \verb"eig" function. The test matrices are formed such that we take a well-conditioned positive definite Hermitian matrix $M$. Then we scale $M$ with an ill-conditioned diagonal matrix $D$ with a decreasing positive diagonal, $A=DMD$. We make sure that $A$ stays Hermitian by setting $A=(A+A^*)/2$ and calculate its exact eigenvalues using the variable-precision arithmetic (vpa) with $100$ significant digits. We show the accuracy of the computed eigenvalues for different block sizes. The results are given for the random pivot strategies $I_{\calO}$, $\calO\in\mathcal{B}_{sp}^{(m)}$, not the same for different block sizes. It can be observed that the complex block Jacobi algorithm is very accurate on this type of test matrices, even when the MATLAB \verb"eig" function is inaccurate.

\section{Conclusion and open problems}

The importance of this work lies in the fact that now we have available tools called the real and complex block Jacobi operators associated with the generalized serial strategies. The real block Jacobi operators have been introduced in~\cite{HaBe17}. These tools can be used to prove the global convergence of the block Jacobi methods for different eigenvalue problems. In this paper, we showed how to use them with the complex block Jacobi methods for the Hermitian, $J$-Hermitian, and normal matrices. We are sure they can be successfully used with the complex block HZ and CJ methods.

Finally, we would like to mention that besides the cyclic pivot strategies, there exist dynamic ordering strategies~\cite{YamOksaVajt89,bec+oks+vaj-15,oks+yam+bec+vaj-18,Kudo+Yas+Yam-18,oks+yam+vaj-19} which are generalizations of the classical pivot strategy~\cite{YamLaKu14,oks+yam+vaj-17}. They are suitable for the block Jacobi methods on parallel machines. It is still an open question how these parallel pivot strategies compare to the known cyclic parallel strategies in terms of efficiency on contemporary parallel computers (cf.~\cite{don+}).

\section*{Acknowledgements}
The authors thank the referees for helpful suggestions which improved the exposition of the paper.


\begin{thebibliography}{10}
\bibitem{bec+oks+vaj-15}
M. Bečka, G. Okša, M. Vajteršic:
\emph{New dynamic orderings for the parallel one-sided block-Jacobi SVD algorithm}.
Parallel Process. Lett. 25(2) (2015) 1550003, 19.
%
\bibitem{BegovicPhD}
E. Begovi\'{c}:
\emph{Convergence of Block Jacobi Methods}.
Ph.D. thesis, University of Zagreb, 2014.
%
\bibitem{BePe23}
E. Begovi\'{c}~Kova\v{c}, A. Perkovi\'{c}:
\emph{Convergence of the Eberlein diagonalization method under the generalized serial pivot strategies}.
Electron. Trans. Numer. Anal. 60 (2024) 238--255.
%
\bibitem{DeVe92}
J. Demmel, K. Veseli\'{c}:
\emph{Jacobi’s method is more accurate than QR}.
SIAM J. Matrix Anal. Appl. 13(4) (1992) 1204--1245.
%
\bibitem{don+}
J. Dongarra, M. Gates, A. Haidar, J. Kurzak, P. Luszczek, S. Tomov, I. Yamazaki:
\emph{The singular value decomposition: anatomy of optimizing an algorithm for extreme scale}.
SIAM Rev. 60(4) (2018) 808--865.
%
\bibitem{DoKoMo09}
Dopico F., Koev P.,Molera J. M.:
\emph{Implicit standard Jacobi gives high relative accuracy}.
Numer. Math. 113 (2009) 519-553.
%
\bibitem{Drmac07}
Z. Drma\v{c}:
\emph{A global convergence proof of cyclic Jacobi methods with block rotations}.
SIAM J. Matrix Anal. Appl. 31(3) (2009) 1329--1350.
%
\bibitem{DrVe08a}
Z. Drma\v{c}, K. Veseli\'{c}:
\emph{New fast and accurate Jacobi SVD algorithm I}.
SIAM J. Matrix Anal. Appl. 29(4) (2008) 1322--1342.
%
\bibitem{DrVe08b}
Z. Drma\v{c}, K. Veseli\'{c}:
\emph{New fast and accurate Jacobi SVD algorithm II}.
SIAM J. Matrix Anal. Appl. 29(4) (2008) 1343--1362.
%
\bibitem{ElIk98}
L. Elsner, K. D. Ikramov:
\emph{Normal matrices: an update}.
Linear Algebra Appl. 285 (1998) 291--303.
%
\bibitem{FoHe60}
G.~E. Forsythe, P. Henrici:
\emph{The cyclic Jacobi method for computing the principal values of a complex matrix}.
Trans. Amer. Math. Soc. 94 (1960) 1--23.
%
\bibitem{GJSW87}
R. Grone, C.R. Johnson, E.M. Sa, H. Wolkowicz:
\emph{Normal matrices}.
Linear Algebra Appl. 87 (1987) 213--225.
%
\bibitem{Han63}
E. R. Hansen:
\emph{On cyclic Jacobi methods}.
SIAM J. Appl. Math. 11 (1963) 449--459.
%
\bibitem{Hari82}
V. Hari:
\emph{On the global convergence of Eberlein method for real matrices}.
Numer. Math. 39 (1982) 361--369.
%
\bibitem{Hari86}
V. Hari:
\emph{On the convergence of cyclic Jacobi-like processes}.
Linear Algebra Appl. 81 (1986) 105--127.
%
\bibitem{Hari91}
V. Hari:
\emph{On pairs of almost diagonal matrices}.
Linear Algebra and Its Appl. 148 (1991) 193--223.
%
\bibitem{Hari91a}
V. Hari: \emph{On sharp quadratic convergence bounds for the serial Jacobi methods}.
Numer. Math. 60 (1991) 375--406.
%
\bibitem{Hari07}
V. Hari:
\emph{Convergence of a block-oriented quasi-cyclic Jacobi method}.
SIAM J. Matrix Anal. Appl. 29(2) (2007) 349--369.
%
\bibitem{Hari09}
V. Hari:
\emph{On block Jacobi annihilators}.
Proceedings of ALGORITMY 2009. Editors: A. Handlovi\v{c}ova et. al. Vysok\'{e} Tatry -- Podbansk\'{e}, Slovakia : Slovak University of Technology in Bratislava, Publishing House of STU, 2009, 429--439.
%
\bibitem{Hari15}
V. Hari:
\emph{Convergence to diagonal form of block Jacobi-type methods}.
Numer. Math. 129(3) (2015) 449--481.
%
\bibitem{Hari18}
V. Hari:
\emph{Complex Cholesky-Jacobi algorithm for PGEP}.
In International conference of numerical Analysis and applied mathematics (ICNAAM 2018), vol.~2116, AIP conference proceedings, 2019.
%
\bibitem{Hari19}
V. Hari:
\emph{On the global convergence of the complex HZ method}.
SIAM J. Matrix Anal. Appl. 40(4) (2019) 1291--1310.
%
\bibitem{Hari21}
V. Hari:
\emph{On the global convergence of the block Jacobi method for the positive definite generalized eigenvalue problem}.
Calcolo 58:24 (2021)
%
\bibitem{Hari22}
V. Hari:
\emph{On the quadratic convergence of the complex HZ method for the positive definite generalized eigenvalue problem}.
Linear Algebra Appl. 632(1) (2022) 153--192.
%
\bibitem{HaBe17}
V. Hari, E. Begovi\'{c}~Kova\v{c}:
\emph{Convergence of the cyclic and quasi-cyclic block Jacobi methods}.
Electron. Trans. Numer. Anal. 46 (2017) 107--147.
%
\bibitem{HaBe21}
V. Hari, E. Begovi\'{c}~Kova\v{c}: \emph{On the convergence of complex Jacobi methods}.
Linear Multilinear Algebra 69(3) (2021) 489--514.
%
\bibitem{HaSiSi10}
V. Hari, S. Singer, S. Singer:
\emph{Block-oriented $J$-{J}acobi method for Hermitian matrices}.
Linear Algebra Appl. 433 (2010) 1491--1512.
%
\bibitem{HaSiSi14}
V. Hari, S. Singer, S. Singer:
\emph{Full block $J$-{J}acobi method for Hermitian matrices}.
Linear Algebra Appl. 444 (2014) 1--27.
%
\bibitem{HeZi68}
P. Henrici, K. Zimmermann:
\emph{An estimate for the norms of certain cyclic Jacobi operators}.
Linear Algebra Appl. 1(4) (1968) 489--501.
%
\bibitem{Kudo+Yas+Yam-18}
S. Kudo, K. Yasuda, Y. Yamamoto:
\emph{Performance of the parallel block Jacobi method with dynamic ordering for the symmetric eigenvalue problem}.
JSIAM Lett. 10 (2018) 41--44.
%
\bibitem{LukPark89}
F. T. Luk, H. Park:
\emph{A proof of convergence for two parallel Jacobi SVD algorithms.}
IEEE Transactions on Computers 38(6) (1989) 806-811.
%
\bibitem{Matejas09}
J. Mateja\v{s}:
\emph{Accuracy of the Jacobi method on scaled diagonally dominant symmetric matrices}.
SIAM J. Matrix Anal. Appl. 31(1) (2009) 133--153.
%
\bibitem{Masc95}
W. F. Mascarenhas: \emph{On the convergence of the Jacobi methods for arbitrary orderings}.
SIAM J. Mat. Anal. Appl. 16(4) (1995) 1197--1209.
%
\bibitem{Naz75}
L. Nazareth:
\emph{On the convergence of the cyclic Jacobi methods}.
Linear Algebra Appl. 12(2) (1975) 151--164.
%
\bibitem{oks+yam+vaj-17}
G. Ok\v{s}a, Y. Yamamoto, M. Vajter\v{s}ic:
\emph{Asymptotic quadratic convergence of the serial block-Jacobi EVD algorithm for Hermitian matrices}.
Numer. Math. 136 (2017) 1071--1095.
%
\bibitem{oks+yam+bec+vaj-18}
G. Ok\v{s}a, Y. Yamamoto, M. Be\v{c}ka, M. Vajter\v{s}ic:
\emph{Asymptotic quadratic convergence of the parallel block-Jacobi EVD algorithm with dynamic ordering for Hermitian matrices}.
BIT 58(4) (2018) 1099--1123.
%
\bibitem{oks+yam+vaj-19}
G. Ok\v{s}a, Y. Yamamoto, M. Be\v{c}ka, M. Vajter\v{s}ic:
\emph{Asymptotic quadratic convergence of the two-sided serial and parallel block-Jacobi SVD algorithm}.
SIAM J. Matrix Anal. Appl. 40(2) (2019) 639--671.
%
\bibitem{RheeHa93}
H. N. Rhee, V. Hari:
\emph{On the global and cubic convergence of a quasy-cyclic Jacobi method}.
Numer. Math. 66 (1993) 97--122
%
\bibitem{Sam71}
A. H. Sameh:
\emph{On Jacobi and Jacobi-like algorithms for parallel computer}.
Math. Comp. 25 (1971) 579--590.
%
\bibitem{Slap92}
I. Slapni{\v{c}}ar:
\emph{Accurate symmetric eigenreduction by a {J}acobi method}.
Ph.D. thesis, Fern{U}niversit{\"{a}}t--{G}esamthochschule, Hagen (1992).
%
\bibitem{Slap03}
I. Slapni{\v{c}}ar:
\emph{Highly accurate symmetric eigenvalue decomposition and hyperbolic {SVD}}.
Linear Algebra Appl. 358 (2003) 387--424.
%
\bibitem{SlapTruh03}
I. Slapni{\v{c}}ar, N. Truhar:
\emph{Relative perturbation theory for hyperbolic singular value problem}.
Linear Algebra Appl. 358 (2003) 367--386.
%
\bibitem{Veselic93}
K. Veseli{\'{c}}:
\emph{A {J}acobi eigenreduction algorithm for definite matrix pairs}.
Numer. Math. 64(1) (1993) 241--269.
%
\bibitem{Veselic00}
K. Veseli{\'{c}}:
\emph{Perturbation theory for the eigenvalues of factorised symmetric matrices}.
Linear Algebra Appl. 309 (2000) 85--102.
%
\bibitem{Wil62}
J. H. Wilkinson:
\emph{Note on the quadratic convergence of the cyclic {J}acobi process}.
Numer. Math. 4 (1962) 296--300.
%
\bibitem{YamOksaVajt89}
Y. Yamamoto, G. Ok\v{s}a, M. Vajter\v{s}ic:
\emph{On convergence to eigenvalues and eigenvectors in the block-Jacobi EVD algorithm with dynamic ordering}.
Linear Algebra Appl. 10(3) (1989) 326--346.
%
\bibitem{YamLaKu14}
Y. Yamamoto, Z. Lan, S. Kudo:
\emph{Convergence analysis of the parallel classical block Jacobi method for the symmetric eigenvalue problem}.
JSIAM Lett. 6 (2014) 57-60.
%
\bibitem{FuZe11}
F. Zhang:
\emph{Matrix Theory}.
Springer, New York, 2011. 
\end{thebibliography}
\end{document}